\documentclass[a4paper,12pt]{amsart}

\usepackage{amsmath}
\usepackage{amsthm}
\usepackage{amssymb}
\usepackage{amscd}
\usepackage{mathrsfs}
\usepackage[all]{xy}
\usepackage{enumerate}

\numberwithin{equation}{section}

\theoremstyle{plain}
\newtheorem{thm}{Theorem}[section]
\newtheorem{dfn}[thm]{Definition}
\newtheorem{exa}[thm]{Example}
\newtheorem{prop}[thm]{Proposition}

\newtheorem{lem}[thm]{Lemma}
\newtheorem*{thm*}{Theorem}

\theoremstyle{definition}
\newtheorem{conj}[thm]{Conjecture}
\newtheorem{rem}[thm]{Remark}

\newcommand{\nc}{\newcommand}
\nc{\Prop}{\begin{prop}}
\nc{\enprop}{\end{prop}}

\def\mod{\mathop{\mathrm{mod}}\nolimits}

\def\Ind{\mathop{\mathrm{Ind}}\nolimits}

\def\Ad{\mathop{\mathrm{Ad}}\nolimits}

\def\id{\mathop{\mathrm{id}}\nolimits}
\def\soc{\mathop{\mathrm{soc}}\nolimits}
\def\cosoc{\mathop{\mathrm{cosoc}}\nolimits}
\def\fmod{\mathop{\mathrm{\text{-}\mod}^{\mathrm{fd}}}\nolimits}

\newcommand{\mf}[1]{{\mathfrak{#1}}}
\newcommand{\mb}[1][m]{{\mathbf{#1}}}
\newcommand{\bb}[1]{{\mathbb{#1}}}

\newcommand{\lr}[2]{{\langle #1,#2 \rangle}}

\nc{\on}{\operatorname}
\newcommand{\C}{{\mathbb C}}
\newcommand{\Q}{\mathbb {Q}}

\newcommand{\Z}{{\mathbb Z}}
\newcommand{\B}{{\mathcal{B}}}

\newcommand{\gl}{{\mathfrak{gl}}}

\newcommand{\seteq}{\mathbin{:=}}

\newcommand{\Lemma}{\begin{lem}}
\newcommand{\enlemma}{\end{lem}}

\newcommand{\g}{{\mathfrak{g}}}

\newcommand{\isoto}[1][]%
{{\mathop{\buildrel{\sim}\over\longrightarrow}\limits_{#1}}}
\renewcommand{\hom}{\operatorname{\it \mathscr{H}\kern-.25em om}}

\newcommand{\M}{{\mathcal{M}}}

\newcommand{\eq}{\begin{eqnarray}}
\newcommand{\eneq}{\end{eqnarray}}

\newcommand{\eqn}{\begin{eqnarray*}}
\newcommand{\eneqn}{\end{eqnarray*}}

\newcommand{\QED}{\end{proof}}
\newcommand{\Proof}{\begin{proof}}

\newcommand{\soplus}{\mathop{\mbox{\normalsize$\bigoplus$}}\limits}

\newcommand{\cl}{\colon}

\newcommand{\ba}{\begin{array}}
\newcommand{\ea}{\end{array}}

\newcommand{\epi}{\twoheadrightarrow}

\newcommand{\set}[2]{\left\{#1 \mid #2 \right\}}

\newcommand{\hs}{\hspace*}

\newcommand{\eqsub}{\begin{subequations}\begin{eqnarray}}
\newcommand{\eneqsub}{\end{eqnarray}\end{subequations}}

\newcommand{\ol}{\overline}

\newcommand{\A}{\mathbf{A}}

\renewcommand{\le}{\leqslant}
\renewcommand{\ge}{\geqslant}

\nc{\la}{\lambda}
\nc{\lam}{\lambda}
\nc{\U}[1][\g]{U_q(#1)}
\nc{\te}{\tilde{e}}
\nc{\tei}{\tilde{e}_i}
\nc{\tf}{\tilde{f}}
\nc{\tfi}{\tilde{f}_i}
\nc{\tU}{\widetilde U_q(\g)}
\nc{\tE}{\widetilde{E}}
\nc{\tF}{\widetilde{F}}

\nc{\BZ}{{\mathbb{Z}}}
\nc{\al}{\alpha}
\nc{\qs}{{q}}
\nc{\lan}{\langle}
\nc{\ran}{\rangle}
\nc{\re}{{\mathrm{re}}}
\nc{\wt}{\operatorname{wt}}
\nc{\Uf}[1][\g]{U^-_q(#1)}
\nc{\Ue}{U^+_q(\g)}
\nc{\eps}{\varepsilon}
\nc{\vphi}{\varphi}
\nc{\sphi}{\varphi^*}
\nc{\seps}{\varepsilon^*}

\nc{\nn}{\nonumber}
\def\max{{\mathop{\mathrm{max}}}}
\nc{\vp}{\varpi}

\nc{\cls}{{\operatorname{cl}}}

\nc{\Wt}{{\operatorname{Wt}}}
\nc{\Us}{U'_q(\g)}
\nc{\La}{\Lambda}
\nc{\ro}{{\rm(}}
\nc{\rf}{{\rm)}}
\nc{\norm}{{\mathrm{norm}}}
\nc{\qbox}{\quad\mbox}
\nc{\braid}{{\mathfrak{B}}}
\nc{\dt}[1]{\tilde{\tilde #1}}
\nc{\Sn}{S^{{\mathrm{norm}}}}
\nc{\aff}{{\mathrm{aff}}}
\nc{\rk}{{\mathrm{rk}}}
\nc{\tQ}{\widetilde{Q}}
\nc{\tP}{\widetilde{P}_\theta}
\nc{\tW}{\widetilde{W}}
\nc{\Dyn}{\mathrm{Dyn}}
\nc{\tD}{\widetilde{\Delta}}
\nc{\height}{{\operatorname{ht}}}
\nc{\bl}{\bigl}
\nc{\br}{\bigr}
\nc{\Hecke}{\mathcal{H}}
\nc{\HA}{\Hecke^{\mathrm{A}}}
\nc{\HB}{\Hecke^{\mathrm{B}}}
\nc{\K}{\mathbf{K}}
\newcommand{\scbul}{{\,\raise1pt\hbox{$\scriptscriptstyle\bullet$}\,}}
\nc{\vac}{{\phi}}
\nc{\Bt}[1][\g]{B_\theta(#1)}
\nc{\be}{\begin{enumerate}}
\nc{\ee}{\end{enumerate}}
\nc{\low}{{\mathrm{low}}}
\nc{\upper}{{\mathrm{up}}}
\nc{\lw}{{\mathrm{low}}}
\nc{\Zodd}{\Z_{\mathrm{odd}}}
\nc{\Ft}[1][n]{\mathbb{P}\mathrm{ol}_{#1}}
\nc{\Ftf}[1][n]{\widetilde{\mathbb{P}\mathrm{ol}}_{#1}}
\nc{\KA}{\on{K}^{\mathrm{A}}}
\nc{\KB}{\on{K}^{\mathrm{B}}}
\nc{\Fc}[1][{n,m}]{\mathbf{F}_{#1}}
\nc{\tphi}{\tilde{\varphi}}
\nc{\CO}{\mathscr{O}}
\nc{\pbw}[1]{\lan #1 \ran}
\nc{\dv}[1]{{[#1]}}
\nc{\Mt}{{\M}_\theta}
\nc{\tVt}{\widetilde{V}_\theta(0)}
\nc{\Vt}[1][0]{V_\theta(#1)}
\nc{\Lt}[1][0]{L_\theta(#1)}
\nc{\tvac}{\widetilde{\vac}}
\nc{\ssum}{\mathop{\mbox{\normalsize{${\sum}$}}}\limits}
\nc{\bnum}{\be[{\rm(i)}]}
\nc{\enum}{\ee}
\nc{\Pt}{P_\theta}
\nc{\suml}{\sum\limits}
\nc{\tL}{\widetilde L}
\nc{\ltcr}{\underset{\mathrm{cry}}{<}}
\nc{\lecr}{\underset{\mathrm{cry}}{\le}}
\nc{\tp}{\tilde{\rho}}
\nc{\dpo}{{(\mspace{-3mu}(}}
\nc{\dpf}{{)\mspace{-3mu})}}
\nc{\odd}{\mathrm{odd}}
\nc{\Bz}{B_\theta(0)}
\nc{\Bs}{B_\theta(\la)}

\title[Symmetric crystals]{\textbf{Symmetric crystals \\
and \\
LLTA type conjectures for the affine Hecke algebras of type B}}
\author{Naoya Enomoto}
\address{Research Institute for Mathematical Sciences\\
Kyoto University\\
Kyoto 606--8502, Japan
}
\email{henon@kurims.kyoto-u.ac.jp}

\author{Masaki Kashiwara}
\address{Research Institute for Mathematical Sciences\\
Kyoto University\\
Kyoto 606--8502, Japan
}
\email{masaki@kurims.kyoto-u.ac.jp}
\thanks{The second author is partially supported by 
Grant-in-Aid for Scientific Research (B) 18340007,
Japan Society for the Promotion of Science.}

\begin{document}

\begin{abstract}
In the previous paper \cite{EK},
we formulated a conjecture 
on the relations between certain classes of irreducible representations
of affine Hecke algebras of type B
and symmetric crystals for $\gl_\infty$.
In the first half of this paper (sections 2 and 3), 
we give a survey of the LLTA type theorem of the affine Hecke algebra 
of type $A$. 
In the latter half (sections 4, 5 and 6), 
we review the construction of the symmetric crystals 
and the LLTA type conjectures for the affine Hecke algebra of type $B$.
\end{abstract}
\maketitle

\section{Introduction}
\subsection{}
The Lascoux-Leclerc-Thibon-Ariki theory connects 
the representation theory of the affine Hecke algebra of {\em type $A$}
with representations of the affine quantum enveloping algebra of type $A$.
Recently, we presented the notion of symmetric crystals 
and conjectured that 
certain classes of irreducible representations of 
the affine Hecke algebras of {\em type B} 
are described by symmetric crystals for $\gl_\infty$ or 
$A_{\ell-1}^{(1)}$ (\cite{EK}). 
In this paper, we review the LLTA-theory 
for the affine Hecke algebra of type $A$, 
the symmetric crystals, and then our conjectures 
for the affine Hecke algebra of type $B$. 
{\em For the sake of simplicity, we restrict ourselves in this note
to the case where the parameters of 
the affine Hecke algebras are not a root of unity.} 

This paper is organized as follows. 
In part I (sections 2 and 3), we review the LLTA-theory 
for the affine Hecke algebras of type $A$. 
In section 2, we recall the representation theory of $U_q(\mf{gl}_{\infty})$, especially the PBW basis, the crystal basis and the global basis. 
In section 3, we recall the representation theory of the affine Hecke algebra of type $A$ and state the LLTA-type theorems. 
In part II (sections 4, 5 and 6), 
we explain symmetric crystals for $\mf{gl}_{\infty}$ 
and the LLTA type conjectures for the affine Hecke algebras of type $B$. 
In section 4, we recall the construction of symmetric crystals 
based on \cite{EK} 
and state the conjecture of existence of the crystal basis 
and the global basis. 
In section 5, we explain 
a combinatorial realization of the symmetric crystals 
for $\mf{gl}_{\infty}$ by using the PBW type basis 
and the $\theta$-restricted multisegments. 
This section is a new additional part to the announcement \cite{EK}. 
The details will appear in \cite{EK2}. 
In section 6, we explain 
the representation theory of the affine Hecke algebra of type $B$ 
and state our LLTA-type conjectures for the affine Hecke algebra of type $B$. 
We add proofs of lemmas and propositions in \cite[section 3.4]{EK}.

\subsection{} Let us recall the LLTA-theory for the affine 
Hecke algebra of type $A$. 

The representation theory of quantum enveloping algebras and the representation theory of affine Hecke algebras have developed independently. 
G. Lusztig \cite{L} constructed 
the PBW type basis and canonical basis of $U_q^{-}(\mf{g})$ 
for the $A$, $D$, $E$ cases. 
The second author \cite{K} defined the crystal basis $B(\infty)$ 
and the (lower and upper) global bases 
$\{G^{\text{low}}(b)\}_{b \in B(\infty)}$, 
$\{G^{\text{up}}(b)\}_{b \in B(\infty)}$ of $U_q^{-}(\mf{g})$. 
The lower global basis coincides with Lusztig's canonical basis. 
On the other hand, A.\ V.\ Zelevinsky \cite{Z} gave 
a parametrization of the irreducible representations 
of the affine Hecke algebra of type $A$ by using multisegments. 
Chriss-Ginzburg \cite{CG} and Kazhdan-Lusztig \cite{KL}  constructed all the irreducible representations of the affine Hecke algebras in geometric methods. 

Lascoux-Leclerc-Thibon conjectured in \cite{LLT} 
that certain composition multiplicities 
(called the decomposition numbers) of the Hecke algebra of type $A$ 
can be written by the transition matrices (specialized at $q=1$) 
between the upper global basis and a standard basis of 
the level $1$ fundamental representation 
of $U_q(\widehat{\mf{sl}_\ell})$. 
In \cite{A}, S.\ Ariki generalized and solved the conjecture for the cyclotomic Hecke algebra and the affine Hecke algebra of type $A$ by using the geometric representation theory of the affine Hecke algebra of type $A$.  In \cite{GV}, 
I.\ Grojnowski and M.\ Vazirani proved the multiplicity-one results for the socle of certain restriction functors and the cosocle of certain induction functors on the category of the finite-dimensional representations of the affine Hecke algebras $\mathcal{H}^A$ of type $A$. 
By using these functors, Grojnowski (\cite{G}) gave the crystal 
structure on the set of irreducible modules over the affine Hecke algebras $\mathcal{H}^A$ of type $A$. In \cite{V}, Vazirani combinatorially constructed the crystal operators on the set of multisegments and proved the compatibility between her actions and Grojnowski's actions. 

For $p \in \bb{C}^*$, let $\mathcal{H}^A_n(p)$ be the affine Hecke algebra of type A of degree $n$ generated by $T_i \ (1 \le i \le n-1)$ and $X_j^{\pm{1}} \ (1 \le j \le n)$. For a subset $J$ of $\bb{C}^*$, we say that a finite-dimensional $\mathcal{H}_n^A$-module is of type $J$ if all the eigenvalues of $X_j \ (1 \le j \le n)$ belong to $J$. We can prove that in order to study the irreducible modules over the affine Hecke algebras of type A, it is enough to treat 
those of type $J$ for an orbit $J$ with respect to the $\Z$-action on $\C^*$ generated by $a \mapsto ap^2$ 
(see Lemma \ref{lem:typeJ}).  
For a $\Z$-orbit $J$, let $K_J(\mathcal{H}_n^A)$ be 
the Grothendieck group of the abelian category of finite-dimensional 
$\mathcal{H}^A_n$-modules of type J, 
and $K_J^A=\oplus_{n\ge0}K_J(\mathcal{H}_n^A)$. 
The LLTA-theory gives the following correspondence 
between the notions in the representation theory 
of a quantum enveloping algebra $U_q(\mf{gl}_{\infty})$ 
and the ones in the representation theory of affine Hecke algebras
of type $A$. 
\renewcommand{\arraystretch}{1.27}
\begin{figure}[h]
\[
\begin{array}{|c|c|}
\hline
\text{the quantum enveloping algebra} & \text{the affine Hecke algebra of type} \ A \\
U_q(\mf{gl}_{\infty}) & \mathcal{H}_n^A(p) \ (n \ge 0) \\
\hline\hline
U_q^{-}(\mf{gl}_{\infty}) & K_J^A=\oplus_{n \ge 0}K_{J}(\mathcal{H}_n^A(p)) \\
\hline
e_a',f_a & \text{certain restrictions $e_a$ and inductions $f_a$} \\
\hline
\text{the crystal basis} \ B(\infty) & \mathcal{M}=\{\text{the multisegments}\} \\
\hline
\text{the upper global basis} & \text{the irreducible modules} \\
\{G^{\text{up}}(b)\}_{b \in B(\infty)} & \{L_b\}_{b \in B(\infty)} \\
\hline
\text{the modified root operators} & \widetilde{e}_a=\soc(e_a),\widetilde{f}_a=\cosoc(f_a) \\
\widetilde{e}_a,\widetilde{f}_a & \widetilde{e}_aL_b=L_{\widetilde{e}_ab},\widetilde{f}_aL_b=L_{\widetilde{f}_ab} \\
\hline
\text{the PBW basis $\{P(b)\}_{b \in B(\infty)}$} 
& \text{the standard modules $\{M(b)\}_{b \in B(\infty)}$} \\
\hline
\end{array}
\]
\caption{Lascoux-Leclerc-Thibon-Ariki correspondence in type A}
\end{figure} 
\renewcommand{\arraystretch}{1.0}

The additive group $\KA_J$ has a structure of Hopf algebra by the restriction and
the induction.
The set $J$ may be regarded as a Dynkin diagram
with $J$ as the set of vertices and with edges between $a\in J$ and $ap^2$.
Let $\g_J$ be the associated Lie algebra, and
$\g_J^-$ the unipotent Lie subalgebra. 
Hence $\g_J$ is isomorphic to $\gl_\infty$ if $p$ has an infinite order.
Let $U_J$ be the group associated
to $\g_J^-$.
Then $\C \otimes \KA_J$ is isomorphic to the algebra $\CO(U_J)$
of regular functions on $U_J$.
Let $\U[\g_J]$ be the associated quantized enveloping algebra.
Then $\Uf[\g_J]$ has a crystal basis $B(\infty)$ and an upper global basis 
$\{G^{\text{up}}(b)\}_{b\in B(\infty)}$.
By specializing $\soplus \C[q,q^{-1}]G^{\text{up}}(b)$
at $q=1$, we obtain $\CO(U_J)$.
Then the LLTA-theory says that the
elements associated to the irreducible $\mathcal{H}^A$-modules
correspond to the image of the upper global basis. Namely, each $b \in B(\infty)$, an irreducible $\mathcal{H}^A$-module $L_b$ is associated and we have 
\[
[e_aL_b:L_{b'}]=e_{a,b,b'}'\vert_{q=1}, 
\quad [f_aL_b:L_{b'}]=f_{a,b,b'}\vert_{q=1}.
\]
Here $[e_aL_b:L_{b'}]$ and $[f_aL_b:L_{b'}]$ are the composition multiplicities of $L_{b'}$ of $e_aL_b$ and $f_aL_b$ in $K_J^A$. (For the definition of the functors $e_a$ and $f_a$ for $a \in J$, see Definition \ref{def:Aef}.) 
The Laurent polynomials $e_{a,b,b'}'$ and $f_{a,b,b'}$ are defined by 
\begin{eqnarray*}
e_a'G^{\text{up}}(b)=\sum_{b' \in B(\infty)}e'_{a,b,b'}G^{\text{up}}(b'), \quad 
f_aG^{\text{up}}(b)=\sum_{b' \in B(\infty)}f_{a,b,b'}G^{\text{up}}(b'). 
\end{eqnarray*}

\subsection{}
Let us explain our analogous conjectures 
for the affine Hecke algebras of type $B$. 

For $p_0, p_1 \in \bb{C}^*$, let $\mathcal{H}_n^B(p_0,p_1)$ be the affine Hecke algebra of type $B$ generated by $T_i \ (0 \le i \le n-1)$ and $X_j \ (1 \le j \le n)$. 
The representation theory of $\mathcal{H}_n^B(p_0,p_1)$ of type $B$ are studied by V.\ Miemietz and Syu Kato. 
In \cite{M}, V.\ Miemietz defined 
certain restriction functors $E_a$ and the induction functors $F_a$ 
on the category of the finite-dimensional representations 
of the affine Hecke algebras of type $B$, 
which are analogous to Grojnowski-Vazirani's construction, 
and proved the multiplicity-one results 
(see sections 6.3 and 6.4). On the other hand, 
S.\ Kato obtained in \cite{Kat} 
a geometric parametrization of the irreducible representations of the affine Hecke algebra $\mathcal{H}_n^B(p_0,p_1)$, 
which is an analogue to 
geometric methods of Kazhdan-Lusztig and Chriss-Ginzburg. 

We say that a finite-dimensional $\mathcal{H}_n^B$-module is of 
type $J\subset \C^*$ 
if all the eigenvalues of $X_j \ (1 \le j \le n)$ belong to $J$.  
Let us consider the $\Z\rtimes \bb{Z}_2$-action on $\bb{C}^*$ 
generated by $a \mapsto ap_1^2$ and $a \mapsto a^{-1}$. 
We can prove that in order to study $\mathcal{H}^B$-modules, it is enough to
study irreducible modules of type $J$ for a $\Z\rtimes\Z_2$-orbit
$J$ in $\C^*$ such that $J$ is a $\Z$-orbit or $J$
contains one of $\pm1,\pm p_0$ (see Proposition \ref{prop:bkB}). 
Let $I=\Z_{\odd}$ be the set of odd integers. 
In this paper, we consider the case 
$J=\set{p_1^k}{k \in I}$ such that $\pm{1}$, $\pm{p_0} \notin J$. 
Let $K_J(\mathcal{H}_n^B)$ be the Grothendieck group of the abelian category of finite-dimensional representations of $\mathcal{H}_n^B(p_0,p_1)$ of type $J$. 

Let $\al_a$ ($a\in J$) be the simple roots with
\[
(\al_a,\al_b)=\begin{cases}
2&\text{if $a=b$,}\\
-1&\text{ if $b=ap_1^{\pm2}$,}\\
0&\text{otherwise.}
\end{cases}
\]
Then the corresponding Lie algebra is $\gl_\infty$.
Let $\theta$ be the involution of $J$ given by $\theta(a)=a^{-1}$. 
In sections 4 and 5, 
we introduce the ring $\mathcal{B}_\theta(\mf{gl}_{\infty})$ 
and the $\mathcal{B}_\theta(\mf{gl}_{\infty})$-module $V_\theta(0)$. 
They are analogues of the reduced $q$-analogue 
$\mathcal{B}_q(\mf{gl}_\infty)$ generated by $e'_a$ and $f_a$,
and the $\mathcal{B}_q(\mf{gl}_{\infty})$-module $U_q^-(\mf{gl}_\infty)$. 
We can prove that $V_\theta(0)$ has the PBW type basis 
$\{P_\theta(b)\}_{b \in B_\theta(0)}$, 
the crystal basis $(L_\theta(0),B_\theta(0))$, 
the lower global basis $\{G_\theta^\low(b)\}_{b \in \mathbb{B}_\theta(0)}$ 
and the upper global basis $\{G_\theta^\upper(b)\}_{b \in B_\theta(0)}$. 
Moreover we can combinatorially describe 
the crystal structure by using the $\theta$-restricted multisegments. 

We conjecture that the irreducible $\mathcal{H}^B$-modules of type $J$ 
are parametrized by $B_\theta(0)$ 
and if $L_b$ is an irreducible $\mathcal{H}^B$-module associated to 
$b \in B_\theta(0)$, 
then we have $\widetilde{E}_aL_b=L_{\widetilde{E}_ab}$,
$\widetilde{F}_aL_b=L_{\widetilde{F}_ab}$ 
and $[E_aL_b:L_{b'}]=E_{a,b,b'}\vert_{q=1}$, 
$[F_aL_b:L_{b'}]=F_{a,b,b'}\vert_{q=1}$. 
(For the definition of the functors $E_a$, $F_a$, $\widetilde{E}_a$ 
and $\widetilde{F}_a$ for $a \in J$, see Definition \ref{def:Bef}.) 
Here the Laurent polynomials $E_{a,b,b'}$ and $F_{a,b,b'}$ are defined by
\[
E_aG_\theta^\upper(b)=\smash{\sum_{b' \in B_\theta(0)}}
E_{a,b,b'}G_\theta^\upper(b'), \quad
F_aG_\theta^\upper(b)=\smash{\sum_{b' \in B_\theta(0)}}
F_{a,b,b'}G_\theta^\upper(b').
\]
\renewcommand{\arraystretch}{1.35}
\begin{figure}[h]
\[
\begin{array}{|c|c|}
\hline
\text{the quantum enveloping algebra} & \text{the affine Hecke algebra of type} \ B \\
U_q(\mf{gl}_{\infty}) \ \text{with} \ \theta & \mathcal{H}_n^B(p_0,p_1) \ (n \ge 0) \\
\hline\hline
V_\theta(0)=U_q^-(\mf{gl}_{\infty})/\sum_i
{U_q^-(\mf{gl}_{\infty})(f_i-f_{\theta(i)})} & K_J^B=\oplus_{n \ge 0}K_J(\mathcal{H}_n^B(p_0,p_1)) \\
\hline
E_a,F_a & \text{certain inductions $E_a$ and restrictions $F_a$} \\
\hline
\text{the crystal basis} \ B_\theta(0) & \mathcal{M}_\theta=\{\text{the $\theta$-restricted multisegments}\} \\
\hline
\text{the upper global basis
$\{G^{\text{up}}_\theta(b)\}_{b \in B_\theta(0)}$} 
& \text{the irreducible modules $\{L_b\}_{b \in B_\theta(0)}$} \\
\hline
\text{the modified root operators} & \widetilde{E}_a=\soc(E_a),\widetilde{F}_a=\cosoc(F_a) \\
\widetilde{E}_a,\widetilde{F}_a & \widetilde{E}_aL_b=L_{\widetilde{E}_ab},\widetilde{F}_aL_b=L_{\widetilde{F}_ab} \\
\hline
\text{the PBW basis $\{P_\theta(b)\}_{b \in B_\theta(0)}$} 
& \text{the standard modules} \\
\hline
\end{array}
\]
\caption{Conjectural correspondence in type B}
\end{figure}
\renewcommand{\arraystretch}{1.0}

\renewcommand{\thepart}{\Roman{part}}
\part{Review on Lascoux-Leclerc-Thibon-Ariki Theory}

\section{Representation Theory of $U_q(\mf{gl}_{\infty})$}

\subsection{Quantized universal enveloping algebras 
and its reduced $q$-analogues}\label{sec:quant}
We shall recall the quantized universal enveloping algebra
$U_q(\g)$.
Let $I$ be an index set (for simple roots),
and $Q$ the free $\Z$-module with a basis $\{\al_i\}_{i\in I}$.
Let $(\scbul,\scbul)\cl Q\times Q\to\Z$ be
a symmetric bilinear form such that
$(\al_i,\al_i)/2\in\Z_{>0}$ for any $i$ and
$(\al_i^\vee,\al_j)\in\Z_{\le0}$ for $i\not=j$ where
$\al_i^\vee\seteq2\al_i/(\al_i,\al_i)$.
Let $q$ be an indeterminate and set 
$\K\seteq\Q(\qs)$.
We define its subrings $\A_0$, $\A_\infty$ and $\A$ as follows.
\begin{eqnarray*}
\A_0&=&\set{f\in\K}{\text{$f$ is regular at $q=0$}},\\
\A_\infty&=&\set{f\in\K}%
{\text{$f$ is regular at $q=\infty$}},\\
\A&=&\Q[\qs,\qs^{-1}].
\end{eqnarray*}
\begin{dfn}\label{U_q(g)}
The quantized universal enveloping algebra $U_q(\g)$ is the $\K$-algebra
generated by elements $e_i,f_i$ and invertible
elements $t_i\ (i\in I)$
with the following defining relations.
\begin{enumerate}[{\rm(1)}]

\item The $t_i$'s commute with each other.

\item
$t_je_i\,t_j^{-1}=q^{(\al_j,\al_i)}\,e_i\ $
and $\ t_jf_it_j^{-1}=q^{-(\al_j,\al_i)}f_i\ $
for any $i,j\in I$.

\item\label{even} $\lbrack e_i,f_j\rbrack
=\delta_{ij}\dfrac{t_i-t_i^{-1}}{q_i-q_i^{-1}}$
for $i$, $j\in I$. Here $q_i\seteq q^{(\al_i,\al_i)/2}$.

\item {\rm(}{\em Serre relation}{\rm)} For $i\not= j$,
\begin{eqnarray*}
\sum^b_{k=0}(-1)^ke^{(k)}_ie_je^{(b-k)}_i=0, \ 
\sum^b_{k=0}(-1)^kf^{(k)}_i
f_jf_i^{(b-k)}=0.
\end{eqnarray*}
Here $b=1-(\al_i^\vee,\al_j)$ and
\begin{eqnarray*}
\ba{l}
e^{(k)}_i=e^k_i/\lbrack k\rbrack_i!\,,\; f^{(k)}_i=f^k_i/\lbrack k
\rbrack_i!\ , \ \lbrack k\rbrack_i=(q^k_i-q^{-k}_i)/(q_i-q^{-1}_i)\,, \ 
\lbrack k\rbrack_i!=\lbrack 1\rbrack_i\cdots \lbrack k\rbrack_i\,.
\ea
\end{eqnarray*}
\end{enumerate}
\end{dfn}
Let us denote by $\Uf$ 
the subalgebra of $\U$ generated by the $f_i$'s.

Let $e'_i$ and $e^*_i$ be the operators on $\Uf$ defined by
$$[e_i,a]=\dfrac{(e^*_ia)t_i-t_i^{-1}e'_ia}{q_i-q_i^{-1}}\quad(a\in\Uf).$$
These operators satisfy the following formulas similar to derivations:
\eq
&&\ba{l}
e_i'(ab)=e_i'(a)b+(\Ad(t_i)a)e_i'b,\quad\\[2ex]
e_i^*(ab)=ae_i^*b+(e_i^*a)(\Ad(t_i)b).
\ea\label{eq:der}
\eneq
The algebra $\Uf$ has a unique symmetric bilinear form $(\scbul,\scbul)$
such that $(1,1)=1$ and
\[
(e'_ia,b)=(a,f_ib)\quad\text{for any $a,b\in\Uf$.}
\]
It is non-degenerate and satisfies $(e^*_ia,b)=(a,bf_i)$.
Let $\B(\g)$ be the algebra generated by the $e_i'$'s and the $f_i$'s.
The left multiplication of $f_j$, $e'_i$ 
and $e_i^*$ have the commutation relations
\[
e'_if_j=q^{-(\al_i,\al_j)}f_je'_i+\delta_{ij}, \ 
e_i^*f_j=f_je_i^*+\delta_{ij}\Ad(t_i),
\]
and both the $e_i'$'s and the $e^*_i$'s satisfy the Serre relations.

\begin{dfn}
The reduced $q$-analogue $\B(\mf{g})$ of $\mf{g}$ is the $\bb{Q}(q)$-algebra generated by $e_i'$ and $f_i$. 
\end{dfn}

\subsection{Review on crystal bases and global bases}
Since $e_i'$ and $f_i$ satisfy the $q$-boson relation, any element $a \in U_q^{-}(\mf{g})$ can be written uniquely as
\[
a=\smash{\sum_{n \ge 0}}f_i^{(n)}a_n \quad \text{with} \ e_i'a_n=0.
\]
Here $f_i^{(n)}=\dfrac{f_i^n}{[n]_i!}$. 
\begin{dfn}
We define the modified root operators $\widetilde{e}_i$ and $\widetilde{f}_i$ on $U_q^{-}(\mf{g})$ by 
\[
\widetilde{e}_ia=\sum_{n \ge 1}f_i^{(n-1)}a_n, \quad \widetilde{f}_ia=\sum_{n \ge 0}f_i^{(n+1)}a_n.
\] 
\end{dfn}
\begin{thm}[\cite{K}]
We define
\begin{eqnarray*}
L(\infty)&=&\sum_{\ell \ge 0,\;i_1, \ldots ,i_{\ell} \in I}
\A_0\tf_{i_1} \cdots \tf_{i_\ell} \cdot 1
 \subset U_q^{-}(\mf{g}), \\
B(\infty)&=&\set{\tf_{i_1} \cdots \tf_{i_\ell} \cdot 1\; \mod qL(\infty)}
{\ell \ge 0,i_1, \cdots ,i_{\ell} \in I} \subset L(\infty)/qL(\infty).
\end{eqnarray*}
Then we have
\bnum
\item $\widetilde{e}_iL(\infty) \subset L(\infty)$ 
and $\widetilde{f}_{i}L(\infty) \subset L(\infty)$, 
\item $B(\infty)$ is a basis of $L(\infty)/qL(\infty)$,
\item
$\widetilde{f}_iB(\infty) \subset B(\infty)$ 
and $\widetilde{e}_iB(\infty) \subset B(\infty) \cup \{0\}$. 
\enum
We call $(L(\infty),B(\infty))$ the crystal basis of $U_q^{-}(\mf{g})$.
\end{thm}
Let $-$ be the automorphism of $\K$ sending $\qs$ to $\qs^{-1}$.
Then $\ol{\A_0}$ coincides with $\A_\infty$. 

Let $V$ be a vector space over $\K$,
$L_0$ an $A$-submodule of $V$,
$L_\infty$ an $\A_\infty$- submodule, and
$V_\A$ an $\A$-submodule.
Set $E\seteq L_0\cap L_\infty\cap V_\A$.

\begin{dfn}[\cite{K}]
We say that $(L_0,L_\infty,V_\A)$ is {\em balanced}
if each of $L_0$, $L_\infty$ and $V_\A$
generates $V$ as a $\K$-vector space,
and if one of the following equivalent conditions is satisfied.
\bnum
\item
$E \to L_0/\qs L_0$ is an isomorphism,
\item
$E \to L_\infty/\qs^{-1}L_\infty$ is an isomorphism,
\item
$(L_0\cap V_\A)\oplus
(\qs^{-1} L_\infty \cap V_\A) \to V_\A$
      is an isomorphism.
\item
$\A_0\otimes_\Q E \to L_0$, $\A_\infty\otimes_\Q E \to L_\infty$,
        $\A\otimes_\Q E \to V_\A$ and $\K \otimes_\Q E \to V$
are isomorphisms.
\enum
\end{dfn}

Let $-$ be the ring automorphism of $\U$ sending
$\qs$, $t_i$, $e_i$, $f_i$ to $\qs^{-1}$, $t_i^{-1}$, $e_i$, $f_i$.

Let $\U_\A$ be the $\A$-subalgebra of
$\U$ generated by $e_i^{(n)}$, $f_i^{(n)}$
and $t_i$.
Similarly we define
$\Uf_\A$.

\begin{thm}
$(L(\infty),L(\infty)^-,\Uf_\A)$ is balanced.
\end{thm}
Let 
\[
G^{\text{low}}\colon L(\infty)/\qs L(\infty)\isoto 
E\seteq L(\infty)\cap L(\infty)^-
\cap \Uf_\A
\] 
be the inverse of $E\isoto L(\infty)/\qs L(\infty)$.
Then $\set{G^{\text{low}}(b)}{b\in B(\infty)}$ forms a basis of $\Uf$.
We call it a (lower) {\em global basis}.
It is first introduced by G.\ Lusztig (\cite{L})
under the name of ``canonical basis'' for the A, D, E cases.

\begin{dfn}
Let
\[
\set{G^\upper(b)}{b \in B(\infty)}
\]
be the dual basis of $\set{G^{\text{low}}(b)}{b \in B(\infty)}$ 
with respect to the inner product $( \scbul,\scbul)$. 
We call it the upper global basis of $U_q^{-}(\mf{g})$.
\end{dfn}

\subsection{Review on the PBW basis}
In the sequel, we set $I=\bb{Z}_{\text{odd}}$ and 
\[
(\alpha_i,\alpha_j)=\left\{
\begin{array}{cl}
2 & \text{for $i=j$,} \\
-1 & \text{for $j=i \pm 2$,}\\
0 & \text{otherwise,}
\end{array}
\right.
\]
and we consider the corresponding
quantum group $\U[\gl_{\infty}]$. In this case,
we can parametrize the crystal basis $B(\infty)$ by the multisegments. 
We shall recall this parametrization and the PBW basis.

\begin{dfn}
For $i,j\in I$ such that $i\le j$, we define a segment $\pbw{i,j}$ 
as the interval $[i,j] \subset \bb{Z}_{\text{odd}}$. 
A multisegment is a formal finite sum of segments:
\[
\mb=\sum_{i \le j}{m_{ij}}\lr{i}{j}
\]
with $m_{i,j}\in\Z_{\ge0}$.
If $m_{i,j}>0$, we sometimes say that
$\pbw{i,j}$ appears in $\mb$.
We denote sometimes $\pbw{i}$ for $\pbw{i,i}$.
We denote by $\mathcal{M}$ the set of multisegments.
We denote by $\emptyset$ the zero element \ro or the empty multisegment\rf\ 
of $\M$.
\end{dfn}
\begin{dfn}
For two segments $\lr{i_1}{j_1}$ and $\lr{i_2}{j_2}$, we define the ordering $\ge_{\text{PBW}}$ by the following:
\[
\lr{i_1}{j_1} \ge_{\text{PBW}} \lr{i_2}{j_2} \Longleftrightarrow \left\{
\begin{array}{l}
j_1>j_2 \\
\text{or} \\
j_1=j_2 \ \text{and} \ i_1 \ge i_2.
\end{array}
\right.
\]
We call this ordering the PBW ordering. 
\end{dfn}
\begin{exa}
We have $\lr{1}{1}>_{\text{PBW}}\lr{-1}{1}>_{\text{PBW}}\lr{-1}{-1}$.
\end{exa}
\begin{dfn}
We define the element $P(\mb) \in U_q^-(\mf{gl}_{\infty})$ indexed 
by a multisegment $\mb$ as follows:
\be[{\rm(1)}]
\item
for a segment $\lr{i}{j}$, we define the element $\lr{i}{j}\in\Uf[\gl_\infty]$ inductively by
\begin{eqnarray*}
\lr{i}{i}&=&f_i, \\
\lr{i}{j}&=&\lr{i}{j-2}\lr{j}{j}-q\lr{j}{j}\lr{i}{j-2},
\end{eqnarray*}
\item
for a multisegment $\displaystyle \mb=\sum_{i \le j}m_{ij}\lr{i}{j}$, we define
\[
P(\mb)=\mathop{\overrightarrow{\prod}}
\lr{i}{j}^{(m_{ij})}.
\]
Here the product $\overrightarrow{\prod}$ 
is taken over segments appearing in $\mb$ 
from large to small with respect to the PBW ordering.
The element $\lr{i}{j}^{(m_{ij})}$ is the divided power of $\lr{i}{j}$ i.e. 
\[
\lr{i}{j}^{(m_{ij})}=\dfrac{1}{[m_{ij}]!}\lr{i}{j}^{m_{ij}}.
\]
\ee
\end{dfn}
\begin{thm}[\cite{L}]
The set of elements $\set{P(\mb)}{\mb \in \mathcal{M}}$ 
is a basis of the $\mathbf{K}$-vector space $U_q^-(\mf{gl}_{\infty})$. 
Moreover this is a basis of the $\mathbf{A}$-module $U_q^-(\mf{gl}_{\infty})_\A$. 
We call this basis the PBW basis of $U_q^-(\mf{gl}_{\infty})$. 
\end{thm}

\begin{dfn}
For two segments $\lr{i_1}{j_1}$ and $\lr{i_2}{j_2}$, 
we define the ordering $\ge_{\text{cry}}$ by the following:
\[
\lr{i_1}{j_1} \ge_{\text{cry}} \lr{i_2}{j_2} \Leftrightarrow \left\{
\begin{array}{l}
j_1>j_2 \\
\text{or} \\
j_1=j_2 \ \text{and} \ i_1 \le i_2.
\end{array}
\right.
\]
We call this ordering the crystal ordering. For $\mb=\sum_{i\le j}m_{i,j}\pbw{i,j}\in \mathcal{M}$ and
and $\mb'=\sum_{i\le j}m'_{i,j}\pbw{i,j}\in \mathcal{M}$,
we define $\mb'\ltcr\mb$
if there exists a segment $\lr{i_0}{j_0}$ such that $m'_{i_0,j_0}<m_{i_0,j_0}$ and $m'_{i,j}=m_{i,j}$ for any $\lr{i}{j}>_{\text{cry}}\lr{i_0}{j_0}$.
\end{dfn}

\begin{exa}
The crystal ordering is different from the PBW ordering. For example, 
we have $\lr{-1}{1}>_{\text{cry}}\lr{1}{1}>_{\text{cry}}\lr{-1}{-1}$,
while we have $\lr{1}{1}>_{\text{PBW}}\lr{-1}{1}>_{\text{PBW}}\lr{-1}{-1}$.
\end{exa}
\begin{dfn}\label{defKop}
We define the crystal structure on $\mathcal{M}$ as follows:
for $\mb=\sum m_{i,j}\pbw{i,j}\in\mathcal{M}$ and $i\in I$,
set $A_k^{(i)}(\mb)=\sum_{k'\ge k}(m_{i,k'}-m_{i+2,k'+2})$ for $k\ge i$.
Define $\eps_i(\mb)$ as $\max\set{A_k^{(i)}(\mb)}{k\ge i}\ge0$.
\bnum
\item If $\eps_i(\mb)=0$, then define
$\te_i(\mb)=0$.
If $\eps_i(\mb)>0$, let
$k_e$ be the largest $k\ge i$
such that $\eps_i(\mb)=A_k^{(i)}(\mb)$
and define
$\te_i(\mb)=\mb-\pbw{i,k_e}+\delta_{k_e\not=i}\pbw{i+2,k_e}$.
\item
Let $k_f$ be the smallest $k\ge i$
such that $\eps_i(\mb)=A_k^{(i)}(\mb)$
and define
$\tf_i(\mb)=\mb-\delta_{k_f\not=i}\pbw{i+2,k_f}+\pbw{i,k_f}$.
\enum
\end{dfn}
\begin{rem}
For $i \in I$,
the actions of the operators $\widetilde{e}_i$ and $\widetilde{f}_i$ on 
$\mb\in\mathcal{M}$ 
are also described by the following algorithm: 

\be[{Step 1.}]
\item Arrange the segments in $\mb$ in the crystal ordering. 
\item
For each segment $\lr{i}{j}$, write $-$, 
and for each segment $\lr{i+2}{j}$, write $+$. 
\item In the resulting sequence of $+$ and $-$, delete a subsequence of the form $+-$ and keep on deleting until no such subsequence remains. 
\ee
Then we obtain a sequence of the form $-- \cdots -++ \cdots +$. 

\be[{\rm(1)}]
\item
$\varepsilon_i(\mb)$ is the total number of $-$ 
in the resulting sequence. 
\item $\widetilde{f}_i(\mb)$ is given as follows:
\be[{(a)}]
\item
If the leftmost $+$ corresponds to a segment $\lr{i+2}{j}$, 
then replace it with $\lr{i}{j}$.
\item If no $+$ exists, add a segment $\lr{i}{i}$ to $\mb$. 
\ee

\item $\widetilde{e}_i(\mb)$  is given as follows:

\be[{(a)}]
\item
If the rightmost $-$ corresponds to a segment $\lr{i}{j}$, 
then replace it with $\lr{i+2}{j}$. 
\item If no $-$ exists, then $\widetilde{e}_i(\mb)=0$. 
\ee
\ee
\end{rem}

\begin{thm}\label{thmgl}
\bnum
\item
 $\displaystyle L(\infty)=\soplus_{\mb \in \mathcal{M}}\mathbf{A}_0P(\mb)$. 
\item
$B(\infty)
=\set{P(\mb) \mod qL(\infty)}{\mb \in \M}$. 
\item
We have
\begin{eqnarray*}
\widetilde{e}_iP(\mb) &\equiv& 
P(\widetilde{e}_i(\mb)) \quad \mod{qL(\infty)}, \\
\widetilde{f}_iP(\mb) &\equiv& P(\widetilde{f}_i(\mb)) 
\quad \mod{qL(\infty)}.
\end{eqnarray*}
Note that $\widetilde{e}_i$ and $\widetilde{f}_i$ 
in the left-hand-side is the modified root operators.
\item
We have the expansion
\[
\overline{P(\mb)}
\in P(\mb)+\sum_{\mb'\ltcr \mb}\A P(\mb').
\]
\ee
\end{thm}
Therefore we can index the crystal basis by multisegments.
By this theorem we can easily see by a standard argument that
$(L(\infty),L(\infty),\Uf_\A)$ is balanced,
and there exists a unique $G^\lw(\mb)\in L(\infty)\cap\Uf_\A$ 
such that $G^\lw(\mb)^{-}=G^\lw(\mb)$ and 
$G^\lw(\mb)\equiv P(\mb) \mod{qL(\infty)}$.
The basis $\{G^\lw(\mb)\}_{\mb\in\M}$ is a lower global basis.

\section{Representation Theory of $\mathcal{H}_n^A$ and the Lascoux-Leclerc-Thibon-Ariki Theory}
\subsection{The affine Hecke algebra of type $A$}
\begin{dfn}
For $p \in \bb{C}^*$, the affine Hecke algebra $\mathcal{H}_n^A$ of type A is a $\bb{C}$-algebra generated by
\[
T_1, \cdots ,T_{n-1}, X_1^{\pm{1}}, \cdots ,X_n^{\pm{1}}
\]
satisfying the following defining relations:
\be[{\rm(1)}]
\item $X_iX_j=X_jX_i$ for any $1 \le i,j \le n$. 
\item {\rm[The braid relations of type $A$]}
\begin{eqnarray*}
\begin{array}{ll}
T_iT_{i+1}T_i=T_{i+1}T_iT_{i+1} & (1 \le i \le n-2), \\
T_iT_j=T_jT_i & (|i-j|>1).
\end{array}
\end{eqnarray*}
\item {\rm[The Hecke relations]}
\[
(T_i-p)(T_i+p^{-1})=0 \quad (1 \le i \le n-1).
\]
\item {\rm[The Bernstein-Lusztig relations]}
\begin{eqnarray*}
\begin{array}{ll}
T_iX_iT_i=X_{i+1} & (1 \le i \le n-1), \\
T_iX_j=X_jT_i & (j \neq i,i+1).
\end{array}
\end{eqnarray*}
\ee
\end{dfn}
\begin{dfn}
For a finite-dimensional $\mathcal{H}_n^A$-module $M$, let 
\[
M=\bigoplus_{a \in (\bb{C}^*)^n}M_a
\]
be the generalized eigenspace decomposition with respect to $X_1,\ldots,X_n$.
Here
\[
M_a\seteq\set{u \in M}{(X_i-a_i)^Nu=0 \ \text{for any
$1 \le i \le n$ and $N \gg 0$}}
\]
for $a=(a_1, \ldots ,a_n) \in (\bb{C}^*)^n$.
\be[{\rm(1)}]
\item
We say that $M$ is of type J if all the eigenvalues of $X_1, \ldots ,X_n$ 
belong to $J \subset \bb{C}^*$. 
\item
Put
\[
K_J^A\seteq\bigoplus_{n \ge 0}K_{J,n}^A.
\]
Here $K_{J,n}^A$ is the Grothendieck group 
of the abelian category of finite-dimensional $\mathcal{H}_n^A$-modules 
of type $J$.
\item
The group $\bb{Z}$ acts on $\bb{C}^*$ by $\bb{Z} \ni n\cl a \mapsto ap^{2n}$.
\ee
\end{dfn}
\begin{lem}\label{lem:typeJ}
Let $J_1$ and $J_2$ be $\Z$-invariant subsets in $\C^*$ such that
$J_1\cap J_2=\emptyset$.
\begin{enumerate}[{\rm(1)}]
\item
If $M$ is an irreducible $\HA_m$-module of type $J_1$
and $N$ is an irreducible $\HA_n$-module of type $J_2$,
then $\Ind_{\HA_m\otimes\HA_n}^{\HA_{m+n}}(M\otimes N)$
is irreducible of type $J_1 \cup J_2$.
\item
Conversely, if $L$ is an irreducible $\HA_n$-module
of type $J_1\cup J_2$, then there exist $m$ $(0\le m\le n)$,
an irreducible $\HA_m$-module $M$ of type $J_1$
and an irreducible $\HA_{n-m}$-module $N$ of type $J_2$
such that $L$ is isomorphic to $\Ind_{\HA_m\otimes\HA_{n-m}}^{\HA_{n}}(M\otimes N)$.
\end{enumerate}
\end{lem}
Hence in order to study the irreducible modules over the affine Hecke
algebras of type A,
it is enough to treat the irreducible modules of type $J$
for an orbit $J$ with respect to the $\Z$-action on $\C^*$.
 
\subsection{The $a$-restriction and the $a$-induction}
For a $\bb{C}$-algebra $A$, let us denote by $A\fmod$ the abelian category of finite-dimensional $A$-modules.
\begin{dfn}\label{def:Aef}
For $a \in \bb{C}^*$, let us define the functors 
\begin{eqnarray*}
e_a:\mathcal{H}_n^A\fmod \to \mathcal{H}_{n-1}^A\fmod, \quad 
f_a:\mathcal{H}_n^A\fmod \to \mathcal{H}_{n+1}^A\fmod
\end{eqnarray*}
by\/{\rm:}
%
$e_aM$  is the generalized $a$-eigenspace of $M$
with respect to the action of $X_n$, and 
$$f_aM\seteq\Ind_{\mathcal{H}_n^A \otimes 
\bb{C}[X_{n+1}^{\pm{1}}]}^{\mathcal{H}_{n+1}^{A}}M \boxtimes 
\langle a \rangle,$$
where $\langle a \rangle$ is the $1$-dimensional 
representation of $\bb{C}[X_{n+1}^{\pm{1}}]$ defined by $X_{n+1} \mapsto a$.

Moreover, put
\[
\widetilde{e}_aM\seteq\soc{e_aM}, \quad \widetilde{f}_aM\seteq\cosoc{f_aM}
\]
for $a \in \bb{C}^{*}$. Here the socle is the maximal semisimple submodule and the cosocle is the maximal semisimple quotient module.
\end{dfn}
\begin{thm}[Vazirani \cite{V}]
Suppose $M$ is irreducible. Then $\widetilde{f}_aM$ is irreducible, and $\widetilde{e}_aM$ is irreducible or $0$ for any $a \in \bb{C}^*$.
\end{thm}

\subsection{LLTA type theorems for the affine Hecke algebra of type $A$}
In this subsection, we consider the case
\[
J=\set{p^k}{k \in \bb{Z}_{\text{odd}}},
\]
and suppose $p$ is not a root of unity. For short, we shall write $e_i,\widetilde{e}_i,f_i$ and $\widetilde{f}_i$ for $e_{p^i},\widetilde{e}_{p^i},f_{p^i}$ and $\widetilde{f}_{p^i}$, respectively. 

The LLTA type theorem for the affine Hecke algebra of type $A$ 
consists of two parts. 
First is a labeling of finite-dimensional irreducible 
$\mathcal{H}^A$-modules by the crystal $B(\infty)$. 
Second is a description of some composition multiplicities 
by using the upper global basis. 

\begin{thm}[Grojnowski-Vazirani \cite{GV}]
There are complete representatives 
\[
\set{L_b}{b \in B(\infty)}
\]
of the finite-dimensional irreducible $\mathcal{H}^A$-modules of type $J$ 
such that 
\[
\widetilde{e}_{i}L_b=L_{\widetilde{e}_ib}, \quad \widetilde{f}_{i}L_b=L_{\widetilde{f}_ib}
\]
for any $i \in I$. 
\end{thm}

\begin{thm}[Ariki \cite{A}]
For $i \in I=\bb{Z}_{\text{odd}}$, let us define $e'_{i,b,b'}, f_{i,b,b'} \in \bb{C}[q,q^{-1}]$ by the coefficients of the expansions:
\begin{eqnarray*}
e_i'G^{\text{up}}(b)=\sum_{b' \in B(\infty)}e'_{i,b,b'}G^{\text{up}}(b'), \quad 
f_iG^{\text{up}}(b)=\sum_{b' \in B(\infty)}f_{i,b,b'}G^{\text{up}}(b'). 
\end{eqnarray*}
Then
\[
[e_{i}L_b:L_{b'}]=e'_{i,b,b'}\vert_{q=1}, 
\quad [f_{i}L_b:L_{b'}]=f_{i,b,b'}\vert_{q=1}.
\]
Here $[M:N]$ is the composition multiplicity of $N$ in $M$ on $K_J^A$.
\end{thm}

\part{The Symmetric Crystals and some LLTA Type Conjectures for Affine Hecke Algebra of Type $B$}

\section{General Definitions and Conjectures for Symmetric Crystals}
We follow the notations in subsection \ref{sec:quant}.
Let $\theta$ be an automorphism of
$I$ such that $\theta^2=\id$ and 
$(\al_{\theta(i)},\al_{\theta(j)})=(\al_i,\al_j)$.
Hence it extends to an automorphism of the root lattice $Q$
by $\theta(\al_i)=\al_{\theta(i)}$,
and induces an automorphism of $\U$.

\begin{dfn}\label{def:Bt}
Let $\B_\theta(\g)$ be the $\K$-algebra
generated by $E_i$, $F_i$, and
invertible elements $T_i$ \ro$i\in I$\rf\ 
satisfying the following defining relations:
\begin{enumerate}[{\rm(i)}]
\item the $T_i$'s commute with each other,
\item
$T_{\theta(i)}=T_i$ for any $i$,
\item
$T_iE_jT_i^{-1}=q^{(\al_i+\al_{\theta(i)},\al_j)}E_j$ and
$T_iF_jT_i^{-1}=q^{(\al_i+\al_{\theta(i)},-\al_j)}F_j$
for $i,j\in I$,
\item
$E_iF_j=q^{-(\al_i,\al_j)}F_jE_i+
(\delta_{i,j}+\delta_{\theta(i),j}T_i)$
for $i,j\in I$,
\item
the $E_i$'s and the $F_i$'s satisfy the Serre relations.
\end{enumerate}
\end{dfn}

We set $E_i^{(n)}=E_i^n/[n]_i!$
and $F_i^{(n)}=F_i^n/[n]_i!$.

\begin{prop}\label{prop:Vtheta}
\begin{enumerate}[{\rm(i)}]
\item
There exists a $\B_\theta(\g)$-module $V_\theta(\la)$
generated by a non-zero vector $\vac_\la$ such that
\be[{\rm(a)}]
\item
$E_i\vac_\la=0$ for any $i\in I$,
\item
$T_i\vac_\la=q^{(\al_i,\la)}\vac_\la$ for any $i\in I$,
\item
$\set{u\in V_\theta(\la)}{\text{$E_iu=0$ for any $i\in I$}}
=\K\vac_\la$.
\ee
Moreover such a $V_\theta(\la)$ is irreducible and
unique up to an isomorphism.
\item
there exists a unique symmetric bilinear form $(\scbul,\scbul)$
on $V_\theta(\la)$ such that $(\vac_\la,\vac_\la)=1$ and
$(E_iu,v)=(u,F_iv)$ for any $i\in I$ and $u,v\in V_\theta(\la)$,
and it is non-degenerate.
\item There exists an endomorphism $-$ of $\Vt[\la]$
such that $\ol{\vac_\la}=\vac_\la$ and
$\ol{av}=\bar{a}\bar{v}$, $\ol{F_iv}=F_i\bar{v}$ for any $a\in \K$
and $v\in\Vt[\la]$.
\end{enumerate}
\end{prop}

The pair $(\Bt,V_\theta(\la))$
is an analogue of
$(\B(\mf{g}),U_q^-(\g))$.
Such a $V_\theta(\la)$ is constructed as follows.
Let $\Uf\vac'_\la$ and $\Uf\vac''_\la$ be a
copy of a free $\Uf$-module.
We give the structure of a $\Bt$-module on them
as follows: for any $i\in I$ and $a\in\Uf$
\eq&&
\left\{
\ba{rcl}
T_i(a\vac'_\la)&=&q^{(\al_i,\la)}(\Ad(t_it_{\theta(i)})a)\vac'_\la,\\[2pt]
E_i(a\vac'_\la)&=&
\bigl(e'_ia+q^{(\alpha_i,\la)}\Ad(t_i)(e^*_{\theta(i)}a)\bigr)\vac'_\la,\\[2pt]
F_i(a\vac'_\la)&=&(f_ia)\vac'_\la
\ea\right.
\eneq
and
\eq&&\left\{\ba{rcl}
T_i(a\vac''_\la)&=&q^{(\al_i,\la)}(\Ad(t_it_{\theta(i)})a)\vac''_\la,\\[2pt]
E_i(a\vac''_\la)&=&(e'_ia)\vac''_\la,\\[2pt]
F_i(a\vac''_\la)&=&
\bigl(f_ia+q^{(\alpha_i,\la)}(\Ad(t_i)a)f_{\theta(i)}\bigr)\vac''_\la.
\ea\right. \label{eq:u''}
\eneq
Then there exists a unique $\Bt$-linear morphism $\psi\cl
\Uf\vac'_\la\to\Uf\vac''_\la$
sending $\vac'_\la$ to $\vac''_\la$.
Its image $\psi(\Uf\vac'_\la)$ is $V_\theta(\la)$.

Hereafter we assume further that
\[
\text{\em there is no $i\in I$ such that $\theta(i)=i$.}
\]
We conjecture that $V_\theta(\la)$ has a crystal basis.
This means the following.
Since $E_i$ and $F_i$ satisfy the $q$-boson relation, 
we define the modified root operators:
\[
\tE_i(u)=\sum_{n\ge1}F_i^{(n-1)}u_n\ 
\text{and}\ 
\tF_i(u)=\sum_{n\ge0}F_i^{(n+1)}u_n,
\]
when writing $u=\sum_{n\ge0}F_i^{(n)}u_n$ with $E_iu_n=0$.
Let $L_\theta(\la)$ be the $\A_0$-submodule of
$V_\theta(\la)$ generated by $\tF_{i_1}\cdots\tF_{i_\ell}\vac_\la$
($\ell\ge0$ and $i_1,\ldots,i_\ell\in I$\,),
and let $B_\theta(\la)$ be the subset 
\[
\set{\tF_{i_1}\cdots\tF_{i_\ell}\vac_\la\bmod \qs L_\theta(\la)}%
{\text{$\ell\ge0$, $i_1,\ldots, i_\ell\in I$}}
\]
of $L_\theta(\la)/\qs L_\theta(\la)$.
\begin{conj}\label{conj:crystal}
\begin{enumerate}
\item
$\tF_iL_\theta(\la)\subset L_\theta(\la)$
and $\tE_iL_\theta(\la)\subset L_\theta(\la)$,
\item
$B_\theta(\la)$ is a basis of $L_\theta(\la)/\qs L_\theta(\la)$,
\item
$\tF_iB_\theta(\la)\subset B_\theta(\la)$,
and
$\tE_iB_\theta(\la)\subset B_\theta(\la)\sqcup\{0\}$,
\item
$\tF_i\tE_i(b)=b$ for any $b\in  B_\theta(\la)$ such that $\tE_ib\not=0$,
and $\tE_i\tF_i(b)=b$ for any $b\in  B_\theta(\la)$.
\end{enumerate}
\end{conj}
Moreover we conjecture that
$V_\theta(\la)$ has a global crystal basis.
Namely we have
\begin{conj}\label{conj:bal}
$(L_\theta(\la),L_\theta(\la)^-,\Vt[\la]^\lw_\A)$
is balanced.
Here $\Vt[\la]^\lw_\A\seteq\Uf_\A\vac_\la$.
\end{conj}
The dual version is as follows.
As in \cite{K}, we have
\Lemma
Assume {\rm Conjecture~\ref{conj:crystal}}.
Then we have
\bnum
\item
$\Lt[\la]=\set{v\in\Vt[\la]}{(\Lt[\la],v)\subset\A_0}$,
\item
Let $(\scbul,\scbul)_0$ be the $\C$-valued symmetric bilinear form
on $\Lt[\la]/q\Lt[\la]$ induced by $(\scbul,\scbul)$.
Then $B_\theta(\la)$ is an orthonormal basis with respect to
 $(\scbul,\scbul)_0$.
\enum
\enlemma
Let us denote by $\Vt[\la]^{\upper}_\A$
the dual space $\set{v\in\Vt[\la]}{(\Vt[\la]^{\lw}_\A,v)\in\A}$.
Then Conjecture~\ref{conj:bal} is equivalent to the following conjecture.
\begin{conj}\label{conj:gls}
$(L_\theta(\la),c(L_\theta(\la)),\Vt[\la]^\upper_\A)$
is balanced.
\end{conj}
Here $c$ is a unique endomorphism of $\Vt[\la]$ such that
$c(\vac_\la)=\vac_\la$ and $c(av)=\bar{a}c(v)$,
$c(E_iv)=E_ic(v)$ for any $a\in \mathbf{K}$ and $v\in \Vt[\la]$.
We have $(c(v'),v)=\ol{(v',\bar{v})}$ for any $v,v'\in \Vt[\la]$.
	
Note that
$\Vt[\la]^\upper_\A$ is the largest $\A$-submodule $M$ of $\Vt[\la]$ such that
$M$ is invariant by the $E_i^{(n)}$'s and
$M\cap\K \vac_\la=\A\vac_\la$.

By Conjecture~\ref{conj:gls}, $\Lt[\la]\cap c(\Lt[\la])\cap\Vt^\upper\to\Lt[\la]/q\Lt[\la]$
is an isomorphism. Let $G^\upper$ be its inverse.
Then $\{G^\upper(b)\}_{b\in\Bs}$ is a basis of $\Vt[\la]$, which we call
the {\em upper global basis} of $\Vt[\la]$.
Note that $\{G^\upper(b)\}_{b\in\Bs}$ is the dual basis
to $\{G^\lw(b)\}_{b\in\Bs}$ with respect to the inner product of
$\Vt[\la]$.
\section{Symmetric Crystals for $\mf{gl}_{\infty}$}
In this section, we consider the case $\mf{g}=\mf{gl}_{\infty}$ and the Dynkin involution $\theta$ of $I$ defined by $\theta(i)=-i$ for $i \in I=\bb{Z}_{\text{odd}}$.
{\scriptsize
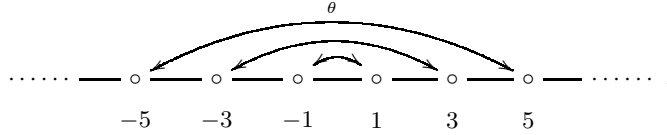
$$
\xymatrix@R=.8ex@C=3ex{
\cdots\cdots\ar@{-}[r]&\circ\ar@{-}[r]
\ar@/^1.8pc/@{<->}[rrrrr]^\theta&\circ\ar@{-}[r]\ar@/^1.2pc/@{<->}[rrr]&\circ\ar@{-}[r]
\ar@/^.7pc/@{<->}[r]&
\circ\ar@{-}[r]&\circ\ar@{-}[r]&\circ\ar@{-}[r]&\cdots\cdots\ .\\
&-5&-3&-1&\;1\;&\;3\;&\;5\;
}
$$}
We shall prove in this case Conjectures \ref{conj:crystal}
and \ref{conj:bal} for $\la=0$.

We set
$$\tVt\seteq\Bt/(\ssum_i\Bt E_i+\ssum_i\Bt(F_i-F_{\theta(i)}))
\simeq
U_q^-(\mf{gl}_{\infty})/\ssum_{i}U_q^-(\mf{gl}_{\infty})(f_i-f_{\theta(i)}).$$
Since $F_i\vac_{0}''=(f_i+f_{\theta(i)})\vac_0''=F_{\theta(i)}\vac_0''$,
we have an epimorphism
\eq
&&\tVt\epi V_\theta(0).
\eneq

It is in fact an isomorphism
(see Theorem~\ref{main:cr}).

\subsection{$\theta$-restricted multisegments}
\begin{dfn}
If a multisegment $\mb$ has the form
\[
\mb=\sum_{-j \le i \le j}m_{ij}\lr{i}{j},
\]
we call $\mb$ a {\em $\theta$-restricted} multisegment. 
We denote by $\mathcal{M}_\theta$ the set of $\theta$-restricted multisegments.
\end{dfn}

\begin{dfn}
 For a $\theta$-restricted segment $\pbw{i,j}$,
we define its modified divided power by
\[
{\lr{i}{j}}^{\dv{m}}=\left\{
\begin{array}{ll}
\pbw{i,j}^{(m)}=\dfrac{1}{[m]!}\lr{i}{j}^m & (i \neq -j), \\
\dfrac{1}{\prod_{\nu=1}^{m}[2\nu]}\lr{-j}{j}^m & (i=-j).
\end{array}
\right.
\]
\end{dfn}

\begin{dfn}
For $\mb \in \mathcal{M}_\theta$, 
we define the elements $P_\theta(\mb) \in \Uf\subset\Bt$ by 
\[
P_\theta(\mb)=\mathop{\overrightarrow{\prod}}\limits
_{\lr{i}{j} \in \mb}\lr{i}{j}^{\dv{m_{ij}}}.
\]
Here the product $\overrightarrow{\prod}$ is taken 
over the segments
appearing in $\mb$ from large to small
with respect to the PBW-ordering. 
\end{dfn}

\subsection{Crystal structure on $\mathcal{M}_\theta$}
\begin{dfn}\label{def:crMt}
Suppose $k>0$. For a $\theta$-restricted multisegment 
$\mb=\suml_{-j\le i\le j}m_{i,j}\pbw{i,j}$, we set
\[
\varepsilon_{-k}(\mb)=\max\set{A_\ell^{(-k)}(\mb)}
{\ell\ge -k},
\]
where
\begin{eqnarray*}
A_\ell^{(-k)}(\mb)&=&\sum_{\ell'\ge\ell}(m_{-k,\ell}-m_{-k+2,\ell+2})
\quad\text{for $\ell>k$,} \\
A_k^{(-k)}(\mb)&=&\sum_{\ell>k}(m_{-k,\ell}-m_{-k+2,\ell})
+2m_{-k,k}+\delta(\text{$m_{-k+2,k}$ is odd}), \\
A_j^{(-k)}(\mb)&=&\sum_{\ell>k}(m_{-k,\ell}-m_{-k+2,\ell})
+2m_{-k,k}-2m_{-k+2,k-2}+\kern-2ex\sum_{-k+2<i\le j+2}\kern-2ex m_{i,k}
-\kern-2ex\sum_{-k+2<i\le j}\kern-2ex m_{i,k-2} \\
&&\hs{48ex}\text{for $-k+2 \le j\le k-2$.}
\end{eqnarray*}
\bnum
\item
Let $n_f$ be the smallest
$\ell\ge-k+2$, with respect to the ordering
$\cdots> k+2>k>-k+2>\cdots>k-2$, such that $\eps_{-k}(\mb)=A_\ell^{(-k)}(\mb)$.
We define
\eqn
\tF_{-k}(\mb)&=&
\begin{cases}
\mb-\pbw{-k+2,n_f}+\pbw{-k,n_f}&\text{if $n_f>k$,}\\
\mb-\pbw{-k+2,k}+\pbw{-k,k}&\text{if $n_f=k$ and $m_{-k+2,k}$ is odd,}\\
\mb-\delta_{k\not=1}\pbw{-k+2,k-2}+\pbw{-k+2,k}
&\text{if $n_f=k$ and $m_{-k+2,k}$ is even,}\\
\mb-\delta_{n_f\not=k-2}\pbw{n_f+2,k-2}+\pbw{n_f+2,k}
&\text{if $-k+2\le n_f\le k-2$.}
\end{cases}
\eneqn

\item
If $\eps_{-k}(\mb)=0$, then $\tE_{-k}(\mb)=0$.
If $\eps_{-k}(\mb)>0$,  then
let $n_e$ be the largest $\ell\ge-k+2$, 
with respect to the above ordering, 
such that $\eps_{-k}(\mb)=A_\ell^{(-k)}(\mb)$.
We define
\eqn
\tE_{-k}(\mb)=
\begin{cases}
\mb-\pbw{-k,n_e}+\pbw{-k+2,n_e}&\text{if $n_e>k$,}\\
\mb-\pbw{-k,k}+\pbw{-k+2,k}
&\text{if $n_e=k$ and $m_{-k+2,k}$ is even,}\\
\mb-\pbw{-k+2,k}+\delta_{k\not=1}\pbw{-k+2,k-2}
&\text{if $n_e=k$ and $m_{-k+2,k}$ is odd,}\\
\mb-\pbw{n_e+2,k}+\delta_{n_e\not=k-2}\pbw{n_e+2,k-2}
&\text{if $-k+2\le n_e\le k-2$.}
\end{cases}&&
\eneqn
\ee
\end{dfn}

\begin{rem}
For $0<k \in I$, 
the actions of $\widetilde{E}_{-k}$ and $\widetilde{F}_{-k}$ on 
$\mb\in\mathcal{M}_\theta$ are described by the following algorithm.
\be[{\rm Step 1.}]
\item
Arrange segments in $\mb$ of the form
$\pbw{-k,j}$ $(j\ge $k), $\pbw{-k+2,j}$ $(j\ge k-2,0)$,
$\pbw{i,k}$ $(-k\le i\le $k), $\pbw{i,k-2}$ $(-k+2\le i\le k-2)$
in the order
\eqn
&&\cdots,\pbw{-k,k+2},\pbw{-k+2,k+2},\;\pbw{-k,k},
\pbw{-k+2,k},\pbw{-k+2,k-2},\\
&&\hs{10ex}\pbw{-k+4,k},\pbw{-k+4,k-2},\cdots,
\pbw{k-2,k},\pbw{k-2,k-2},\pbw{k}.
\eneqn

\item Write signatures for each segment appearing in 
$\mb$ by the following rules.
\bnum
\item If a segment is not $\lr{-k+2}{k}$, then
\begin{itemize}
\item{For $\lr{-k}{k}$, write $--$,}
\item{For $\lr{-k}{j}$ with $j> k$, write $-$,} 
\item{For $\lr{-k+2}{k-2}$ with $k>1$, write $++$,} 
\item{For $\lr{-k+2}{j}$ with $j>k$, write $+$,}
\item{For $\lr{j}{k}$ if $-k<j\le k$, write $-$,} 
\item{For $\lr{j}{k-2}$ if $-k+2<j\le k-2$, write $+$,}
\item{If otherwise, write no signature.} 
\end{itemize}
\item
For segments $m_{-k+2,k}\lr{-k+2}{k}$, 
if $m_{-k+2,k}$ is even, then write no signature, 
and if $m_{-k+2,k}$ is odd, then write a sequence $-+$. 
\ee
\item
In the resulting sequence of $+$ and $-$, 
delete a subsequence of the form $+-$ 
and keep on deleting until no such subsequence remains. 
\ee
Then we obtain a sequence of the form $-- \cdots -++ \cdots +$. 
\be[{\rm(1)}]
\item
$\varepsilon_{-k}(\mb)$ is given
as the total number of $-$ in the resulting sequence. 
\item $\widetilde{F}_{-k}(\mb)$ is given as follows:
\bnum
\item
if the leftmost $+$ corresponds to a segment $\lr{-k+2}{j} \ (j >k)$, 
then replace the segment with $\lr{-k}{j}$,
\item
if the leftmost $+$ corresponds to a segment $\lr{j}{k-2}$, 
then replace the segment with $\lr{j}{k}$,
\item
f the leftmost $+$ corresponds to segment $\lr{-k+2}{k}^{m_{-k+2,k}}$, 
then replace one of the segments with $\lr{-k}{k}$,
\item if no $+$ exists, add a segment $\lr{k}{k}$ to $\mb$. 
\ee
\item
$\widetilde{E}_{-k}(\mb)$ is given as follows: 
\bnum
\item
if the rightmost $-$ corresponds to a segment $\lr{-k}{j}$, 
then replace the segment with $\lr{-k+2}{j}$,
\item
if the rightmost $-$ corresponds to a segment $\lr{j}{k} \ (j \neq -k+2)$, 
then replace the segment with $\lr{j}{k-2}$,
\item
if the rightmost $-$ corresponds to segments $m_{-k+2,k}\lr{-k+2}{k}$, 
then replace one of the segment with $\lr{-k+2}{k-2}$,
\item
if no $-$ exists, then $\widetilde{E}_{-k}(\mb)=0$. 
\ee
\ee
\end{rem}
\begin{dfn}
For $k\in I_{>0}$, we define
$\tF_{k}$, $\tE_k$ and $\eps_k$ by the same rule as in 
{\rm Definition~\ref{defKop}} for
$\tf_k$ and $\te_k$.
\end{dfn}

\begin{thm}\label{th:Excr}
By $\tF_k$, $\tE_k$, $\eps_k$ \ro $k\in I$\rf,
$\Mt$ is a crystal, 
namely, for any $k\in I$, we have 
\bnum
\item
$\tF_{k}\Mt\subset \Mt$ and
$\tE_{k} \Mt\subset \Mt\sqcup\{0\}$,
\item
$\tF_{k}\tE_{k}(\mb)=\mb$ if $\tE_{k}(\mb)\not=0$, 
and $\tE_{k}\circ\tF_{k}=\id$,
\item 
$\eps_{k}(\mb)=\max\set{n\ge0}{\tE^n(\mb)\not=0}<\infty$ for any 
$\mb\in\Mt$.
\enum
\end{thm}
\begin{exa}
\begin{enumerate}[{\rm(1)}]
\item We shall write $\{a,b\}$ for $a\lr{-1}{1}+b\langle 1 \rangle$. 
The following diagram is the part of the crystal graph of $\Bz$ that concerns only the $1$-arrows and the $(-1)$-arrows.
$$
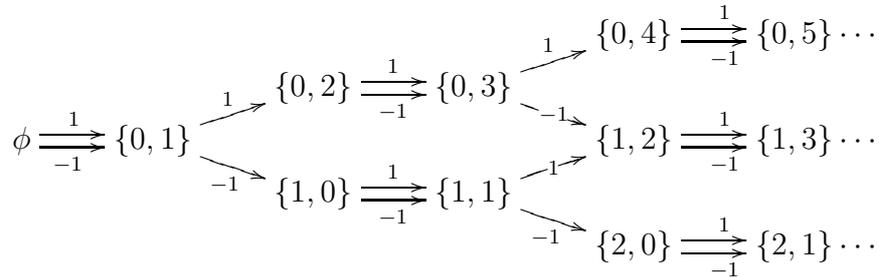
\xymatrix@R=.08em@C=2em{&&&&\{0,4\}\ar@<.2pc>[r]^{1}\ar@<-.2pc>[r]_{-1}&
\{0,5\}\cdots\\
&&\{0,2\}\ar@<.2pc>[r]^{1}\ar@<-.2pc>[r]_{-1}&\{0,3\}
\ar[ru]^{1}\ar[rd]|{-1}\\
\vac\ar@<.2pc>[r]^(.4){1}\ar@<-.2pc>[r]_(.35){-1}&
\{0,1\}\ar[ru]^{1}\ar[rd]_{-1}&&&\{1,2\}\ar@<.2pc>[r]^{1}\ar@<-.2pc>[r]_{-1}&\{1,3\}\cdots\\
&&\{1,0\}\ar@<.2pc>[r]^{1}\ar@<-.2pc>[r]_{-1}&\{1,1\}\ar[ru]|{1}\ar[rd]_{-1}\\
&&&&\{2,0\}\ar@<.2pc>[r]^{1}\ar@<-.2pc>[r]_{-1}&{\{2,1\}\cdots}
}$$
Especially the part of $(-1)$-arrows is the following diagram.
$$\xymatrix@C=6ex{
\{0,2n\}\ar[r]^(.45){-1}&
\{0,2n+1\}\ar[r]^(.55){-1}&
\{1,2n\}\ar[r]^(.45){-1}&
\{1,2n+1\}\ar[r]^(.55){-1}&
\{2,2n\}\ar[r]^(.55){-1}&
\cdots
}$$
\item The following diagram is the part of the crystal graph of 
$B_\theta(0)$ that concerns only the $(-1)$-arrows and the $(-3)$-arrows. 
This diagram is isomorphic as a graph to the crystal graph 
of $A_2$.
$$
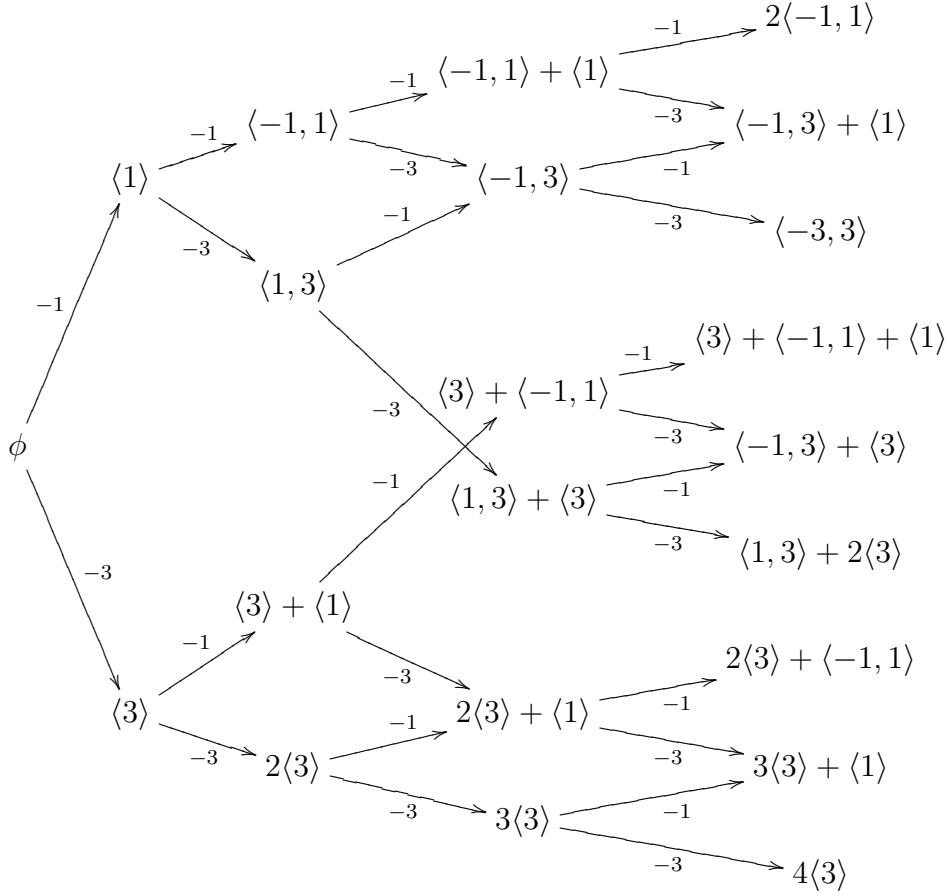
\xymatrix@R=.1em@C=2em{&&&&2\lr{-1}{1}\\
&&&\lr{-1}{1}+\langle 1 \rangle \ar[ru]^{-1}\ar[rd]_{-3}\\
&& \lr{-1}{1}\ar[ru]^{-1}\ar[rd]_{-3} && \lr{-1}{3}+\langle 1 \rangle \\
& \langle 1 \rangle\ar[ru]^{-1}\ar[rdd]_{-3} 
&& \lr{-1}{3}\ar[ru]_{-1}\ar[rd]_{-3} & \\
&&&& \lr{-3}{3} \\
&& \lr{1}{3}\ar[ruu]^{-1}\ar[rdddd]_{-3} & \\
&&&& \pbw{3}+\pbw{-1,1}+\pbw{1} \\
&&& \langle 3 \rangle+\lr{-1}{1}\ar[ru]^(.4){-1}\ar[rd]_{-3} & \\
\vac\ar[ruuuuu]^{-1}\ar[rddddd]^{-3} &&&& \lr{-1}{3}+\langle 3 \rangle \\
&&& \pbw{1,3}+\langle 3 \rangle\ar[ru]_{-1}\ar[rd]_{-3} & \\
&&&& \pbw{1,3}+2\pbw{3}\\
&& \langle 3 \rangle+\langle 1 \rangle\ar[ruuuu]^{-1}\ar[rdd]_{-3} & \\
&&&& 2\langle 3\rangle+\pbw{-1,1}\\
& \langle 3 \rangle\ar[ruu]^{-1}\ar[rd]_{-3} 
&& 2\langle 3 \rangle+\langle 1 \rangle\ar[ru]_{-1}\ar[rd]_{-3} & \\
&& 2\langle 3 \rangle \ar[ru]^{-1}\ar[rd]_{-3} && 3\pbw{3}+\pbw{1}\\
&&&3\langle 3 \rangle \ar[ru]_{-1}\ar[rd]_{-3}\\
&&&&4\langle 3 \rangle
}$$
\item Here is the part of the crystal graph of $\Bz$
that concerns only the $n$-arrows and the $(-n)$-arrows
for an odd integer $n\ge 3$:
$$\xymatrix@C=6ex{
\vac\ar@<.2pc>[r]^(.5){n}\ar@<-.2pc>[r]_(.45){-n}&
\pbw{n}\ar@<.2pc>[r]^(.45){n}\ar@<-.2pc>[r]_(.4){-n}&
2\pbw{n}\ar@<.2pc>[r]^(.45){n}\ar@<-.2pc>[r]_(.4){-n}&
3\pbw{n}\ar@<.2pc>[r]^(.4){n}\ar@<-.2pc>[r]_(.35){-n}&
4\pbw{n}\cdots
}$$
\end{enumerate}
\end{exa}

\subsection{Main Theorem}
We write $\vac$ for the generator $\vac_0$ of $\Vt$, for short.
\begin{thm}\label{main:cr}
\bnum
\item The morphism
$$\tVt=\Uf/\sum_{k\in I}\Uf(f_k-f_{-k})\to \Vt$$
is an isomorphism.
\item
$\{\Pt(\mb)\vac\}_{\mb\in\Mt}$ is a basis of
the $\K$-vector space $\Vt$.
\item
Set 
\eqn
&&\Lt\seteq\sum_{\ell\ge0, i_1,\ldots,i_\ell\in I}
\A_0\tF_{i_1}\cdots\tF_{i_\ell}\vac\subset\Vt,\\
&&\Bz=\set{\tF_{i_1}\cdots\tF_{i_\ell}\vac\mod q\Lt}%
{\ell\ge0, i_1,\ldots,i_\ell\in I}.
\eneqn
Then, $\Bz$ is a basis of $\Lt/q\Lt$ and
$(\Lt,\Bz)$ is a crystal basis of $\Vt$,
and the crystal structure coincide with the one of $\Mt$.
\item 
More precisely, we have
\be[{\rm(a)}]
\item
$\Lt=\sum_{\mb\in\Mt}\A_0\Pt(\mb)\vac$,
\item
$\Bz=\set{\Pt(\mb)\vac\mod q\Lt}{\mb\in\Mt}$,
\item
for any $k\in I$ and $\mb\in\Mt$, we have
\be[{\rm (1)}]
\item
$\tF_k\Pt(\mb)\vac\equiv\Pt(\mb)\vac\;\mod q\Lt$,
\item
$\tE_k\Pt(\mb)\vac\equiv\Pt(\tE_k\mb)\vac\;\mod q\Lt$, where
we understand $\Pt(0)=0$,
\item
$\tE_k^n\Pt(\mb)\vac\in q\Lt$ if and only if $n>\eps_k(\mb)$.
\ee
\ee
\enum
\end{thm}

\subsubsection{Global Basis of $V_\theta(0)$}
Recall that $\A=\Q[q,q^{-1}]$, 
and $\Vt_\A=\Uf[\gl_\infty]_\A\vac$.
\begin{lem}
\bnum
\item $\Vt_\A=\soplus_{\mb\in\Mt}\A \Pt(\mb)\vac$.
\item For $\mb\in\M$,
$$\ol{\Pt(\mb)\vac}\in \Pt(\mb)\vac+\sum_{\mb[n]\ltcr\mb}\A \Pt(\mb[n])\vac.$$
\end{enumerate}
\end{lem}
By the above lemma, we obtain the following theorem.
\begin{thm}
\bnum
\item $(\Lt,\Lt^-,\Vt_\A)$ is balanced.
\item For any $\mb\in\Mt$, there exists a unique $G_\theta^\lw(\mb)
\in \Lt\cap\Vt_\A$
such that
$\ol{G_\theta^\lw(\mb)}=G_\theta^\lw(\mb)$ and $G_\theta^\lw(\mb)\equiv \Pt(\mb)\vac\;\mod q\Lt$.
\item
$G_\theta^\lw(\mb)\in \Pt(\mb)\vac
+\sum_{\mb[n]\ltcr\mb}q\,\C[q]\Pt(\mb[n])\vac$
for any $\mb\in\Mt$.
\end{enumerate}
\end{thm}

\section{Representation Theory of $\mathcal{H}_n^B$ and LLTA Type Conjectures}

\subsection{The affine Hecke algebra of type $B$}
\begin{dfn}
For $p_0$, $p_1 \in \bb{C}^*$, 
the affine Hecke algebra $\mathcal{H}_n^B$ of type B 
is a $\bb{C}$-algebra generated by
\[
T_0,T_1, \cdots ,T_{n-1}, X_1^{\pm{1}}, \cdots ,X_n^{\pm{1}}
\]
satisfying the following defining relations\/{\rm:}
\bnum
\item $X_iX_j=X_jX_i$ for any $1 \le i,j \le n$. 
\item {\rm[The braid relations of type $B$]}
\begin{eqnarray*}
\begin{array}{ll}
T_0T_1T_0T_1=T_1T_0T_1T_0, & \\
T_iT_{i+1}T_i=T_{i+1}T_iT_{i+1} & (1 \le i \le n-2), \\
T_iT_j=T_jT_i & (|i-j|>1).
\end{array}
\end{eqnarray*}
\item {\rm[The Hecke relations]}
\[
(T_0-p_0)(T_0+p_0^{-1})=0, \quad (T_i-p_1)(T_i+p_1^{-1})=0 \quad (1 \le i \le n-1).
\]
\item {\rm[The Bernstein-Lusztig relations]}
\begin{eqnarray*}
\begin{array}{ll}
T_0X_1^{-1}T_0=X_1, & \\
T_iX_iT_i=X_{i+1} & (1 \le i \le n-1), \\
T_iX_j=X_jT_i & (j \neq i,i+1).
\end{array}
\end{eqnarray*}
\ee
\end{dfn}
Note that the subalgebra generated by $T_i \ (1 \le i \le n-1)$ and $X_j^{\pm{1}} \ (1 \le j \le n)$ is isomorphic to the affine Hecke algebra $\mathcal{H}_n^{A}$. 

We assume that $p_0$, $p_1\in\C^*$ satisfy
\[
p_0^2\not=1,\ p_1^2\not=1.
\]
Let us denote by $\Ft$ the Laurent polynomial ring
$\C[X_1^{\pm1},\ldots,X_n^{\pm1}]$,
and by $\Ftf$ its quotient field $\C(X_1,\ldots, X_n)$.
Then $\mathcal{H}_n^B$ is isomorphic to the tensor product
of $\Ft$ and the subalgebra generated by the $T_i$'s
that is isomorphic to the Hecke algebra of type $B_n$.
We have
\[
T_ia=(s_ia)T_i+(p_i-p_i^{-1})\dfrac{a-s_ia}{1-X^{-\al_i^\vee}}
\quad\text{for $a\in\Ft$.}
\]
Here $p_i=p_1$ ($1<i<n$),
and $X^{-\al_i^\vee}=X_1^{-2}$ ($i=0$)
and $X^{-\al_i^\vee}=X_{i}X_{i+1}^{-1}$ ($1\le i<n$).
The $s_i$'s are the Weyl group action on $\Ft$:
$(s_ia)(X_1,\ldots,X_n)=a(X_1^{-1},X_2,\ldots,X_n)$ for $i=0$
and
$(s_ia)(X_1,\ldots,X_n)=a(X_1,\ldots,X_{i+1},X_i,\ldots,X_n)$ for $1\le i<n$.

Note that $\mathcal{H}_n^B=\C$ for $n=0$.

The algebra $\mathcal{H}_n^B$ acts faithfully on
$\mathcal{H}_n^B/\sum_{i=0}^{n-1}\mathcal{H}_n^B(T_i-p_i)\simeq\Ft$.
Set 
\[
\vphi_i=(1-X^{-\al_i^\vee})T_i-(p_i-p_i^{-1}) \in \mathcal{H}_n^B
\]
and 
\[
\tphi_i=(p_i^{-1}-p_iX^{-\al_i^\vee})^{-1}\vphi_i\in\Ftf\otimes_{\Ft}
\mathcal{H}_n^B.
\]
Then the action of $\tphi_i$ on $\Ft$ coincides with $s_i$.
They are called {\em intertwiners}.

\subsection{Block decomposition of $\mathcal{H}_n^B$-mod}
For $n,m\ge 0$, set
$$\Fc\seteq\C[X_1^{\pm1},\ldots,X_{n+m}^{\pm1},D^{-1}],$$
where
\[
D\seteq\prod\limits_{1\le i\le n<j\le n+m}\hs{-2ex}
(X_i-p_1^2 X_j)(X_i-p_1^{-2}X_j)
(X_i-p_1^2 X_j^{-1})(X_i-p_1^{-2}X_j^{-1})(X_i-X_j)(X_i-X_j^{-1}).
\]
Then we can embed
$\mathcal{H}_n^B$ into $\mathcal{H}_{n+m}^B\otimes_{\Ft[n+m]}\Fc$
by 
\[
T_0\mapsto \tphi_n\cdots\tphi_1\,T_0\,\tphi_1\cdots\tphi_n, \quad 
T_i\mapsto T_{i+n}\ (1\le i<m),\quad
X_i\mapsto X_{i+n}\ (1\le i\le m).
\]
Its image commute with $\mathcal{H}_n^{B}\subset\mathcal{H}_{n+m}^B$.
Hence $\mathcal{H}_{n+m}^B\otimes_{\Ft[n+m]}\Fc$
is a right $\mathcal{H}_n^B\otimes\mathcal{H}_m^B$-module. 
Note that $(\mathcal{H}_n^B \otimes \mathcal{H}_m^B) \otimes_{\Ft[n+m]}
\mathbf{F}_{n,m}=\mathbf{F}_{n,m} \otimes_{\Ft[n+m]} 
(\mathcal{H}_n^B \otimes \mathcal{H}_m^B)$ is an algebra. 
\begin{lem}
$\mathcal{H}_{n+m}^A\hs{-1ex}\mathop\otimes\limits_{\mathcal{H}_n^A\otimes
\mathcal{H}_m^A}\hs{-1.5ex}
\bl(\mathcal{H}_n^B\otimes\mathcal{H}_m^B
\br)\otimes_{\Ft[n+m]}{\Fc}\isoto
\mathcal{H}_{n+m}^B\otimes_{\Ft[n+m]}\Fc$.
\end{lem}
\begin{proof}
Let $W_n^A$ and $W_n^B$ be the finite Weyl group of type $A$ and $B$. 
Note that $|W_{n+m}^A| \cdot |W_n^B| \cdot |W_m^B|/(|W_n^A| \cdot |W_m^A|)
=|W_{n+m}^B|$. 
Hence the both sides are free modules of rank 
$|W_{n+m}^B|$ over $\mathbf{F}_{n,m}$. 
We prove that the map is surjective. 

For short, we denote the image of 
$\mathcal{H}_{n+m}^A\hs{-1ex}\mathop\otimes\limits_{\mathcal{H}_n^A
\otimes\mathcal{H}_m^A}\hs{-.5ex}
\bl(\mathcal{H}_n^B\otimes\mathcal{H}_m^B
\br)\otimes_{\Ft[n+m]}{\Fc}$ by $\mathcal{H}_{n,m}^{\text{loc}} \subset \HB_{n+m} \otimes_{\Ft[n+m]} \mathbf{F}_{n,m}$. 
Note that $\tilde{\varphi}_i \cdots \tilde{\varphi}_n \in \HA_{n+m}\otimes_{\Ft[n+m]}\mathbf{F}_{n,m}$ for $1 \le i \le n$.

First, we have $\tphi_n\cdots\tphi_1\,T_0\,\tphi_1\cdots
\tphi_n \in \mathcal{H}_m^B \otimes_{\Ft}\mathbf{F}_{n,m}$. 
Since $(\tphi_n\cdots\tphi_1)^{-1}
=\tphi_1 \cdots \tphi_n 
\in \mathcal{H}_{n+m}^A \otimes_{\Ft}\mathbf{F}_{n,m}$, 
we have $T_0\,\tphi_1\cdots\tphi_n \in \mathcal{H}_{n,m}^{\text{loc}}$. 

Second, note that
\[
T_i=\left(\tphi_i(p_i^{-1}-p_iX_i^{-1}X_{i+1})-(p_i-p_i^{-1})X_i^{-1}X_{i+1}
\right)(1-X_i^{-1}X_{i+1})^{-1} \ (1 \le i<n).
\]
If $T_0T_1 \cdots T_{i-1}\,\tphi_i \cdots \tphi_n 
\in \mathcal{H}_{n,m}^{\text{loc}}$, 
then $T_0T_1 \cdots T_{i}\,\tphi_{i+1} \cdots \tphi_n 
\in \mathcal{H}_{n,m}^{\text{loc}}$ for $1 \le i<n$.
Indeed, we have
\begin{eqnarray*}
T_0 \cdots T_i\,\tphi_{i+1} \cdots \tphi_n
&=&T_0 \cdots T_{i-1}\,\tphi_i \cdots \tphi_n\,
(p_i^{-1}-p_iX_i^{-1}X_{n+1})(1-X_i^{-1}X_{n+1})^{-1} \\
&{}& \quad -(p_i-p_i^{-1})T_0 \cdots T_{i-1}\,\tphi_{i+1} \cdots \tphi_n\,
X_i^{-1}X_{n+1}(1-X_i^{-1}X_{n+1})^{-1}
\end{eqnarray*}
and
\[
T_0 \cdots T_{i-1}\,\tphi_{i+1} \cdots \tphi_n
=\tphi_{i+1} \cdots \tphi_n\,T_0 \cdots T_{i-1} 
\in \HA_{n+m}\mathbf{F}_{n,m}\HB_n.
\]
 Therefore $T_0T_1 \cdots T_n \in \mathcal{H}_{n,m}^{\text{loc}}$. 
Hence $T_0T_1 \cdots T_i \in \mathcal{H}_{n,m}^{\text{loc}}$
($1 \le i<n+m$). 
Indeed, if $i<n$, then $T_0T_1 \cdots T_i \in \mathcal{H}_n^B$. 
If $n \le i$, then $T_0T_1 \cdots T_n \in \mathcal{H}_{n,m}^{\text{loc}}$ 
and $T_{n+1} \cdots T_i \in \mathcal{H}_m^B$.

Finally, we prove the surjectivity by the induction on $m$. Note that 
\[
\mathcal{H}_{n+m}^B=\sum_{i=1}^{n+m}T_iT_{i+1} \cdots T_{n+m-1}\mathcal{H}_{n+m-1}^{B}+\sum_{i=0}^{n+m-1}T_i \cdots T_1T_0T_1 \cdots T_{n+m-1}\HB_{n+m-1}
\]
and $T_iT_{i+1} \cdots T_{n+m-1} \in \mathcal{H}_{n+m-1}^A$. Furthermore, $\mathcal{H}_{n+m-1}^{B} \subset \mathcal{H}_{n,m-1}^{\text{loc}}$ by the induction hypothesis. Thus it is sufficient to prove that $T_0\mathcal{H}_{n+m}^{A,\text{fin}} \subset \mathcal{H}_{n,m}^{\text{loc}}$. Here, $\mathcal{H}_{n+m}^{A,\text{fin}}$ is the subalgebra of $\HA_{n+m}$ generated by $T_1, \ldots ,T_{n+m-1}$. This follows from
\[
\mathcal{H}_{n+m}^{A,\text{fin}}
=\sum_{i=0}^{n+m-1}\langle T_2, \cdots ,T_{n+m-1} 
\rangle T_1T_2 \cdots T_i
\]
and $T_0T_1 \cdots T_i \in \mathcal{H}_{n,m}^{\text{loc}}$. 
\end{proof}
\begin{dfn}
For a finite-dimensional $\mathcal{H}_n^B$-module $M$, let 
\[
M=\bigoplus_{a \in (\bb{C}^*)^n}M_a
\]
be the generalized eigenspace decomposition 
with respect to $X_1, \ldots ,X_n${\rm:}
\[
M_a\seteq\set{u \in M}{(X_i-a_i)^Nu=0 \ 
\text{for any $1 \le i \le n $ and $N \gg 0$}}
\]
for $a=(a_1, \ldots ,a_n) \in (\bb{C}^*)^n$. 
\begin{enumerate}[{\rm(1)}]
\item We say that $M$ is of type $J$ if all the eigenvalues 
of $X_1, \ldots ,X_n$ belong to $J \subset \bb{C}^*$. Put
\[
K_J^B\seteq\bigoplus_{n \ge 0}K_{J,n}^B.
\]
Here $K_{J,n}^B$ is the Grothendieck group of the abelian category of finite-dimensional $\mathcal{H}_n^B$-modules of type $J$.
\item The semi-direct product group
$\Z\rtimes\Z_2=\Z\times\{1,-1\}$ acts on $\C^*$
by $(n,\epsilon)\cl a\mapsto a^\epsilon p_1^{2n}$.
\item Let $J_1$ and $J_2$ be $\Z\rtimes\Z_2$-invariant subsets of
$\C^*$ such that $J_1 \cap J_2=\emptyset$.
Then for an $\mathcal{H}_n^B$-module $N$
of type $J_1$ and an $\mathcal{H}_m^B$-module $M$ of type $J_2$,
the action of $\Ft[n+m]$ on $N\otimes M$ extends to an action of $\Fc$.
We set 
\[
N\diamond M\seteq (\mathcal{H}_{n+m}^B\otimes_{\Ft[n+m]}\Fc)\otimes_{(\mathcal{H}_n^B\otimes\mathcal{H}_m^B)\otimes_{\Ft[n+m]}{\Fc}}(N\otimes M).
\]
\end{enumerate}
\end{dfn}
By the lemma above, $N\diamond M$ is isomorphic to
$\Ind_{\mathcal{H}_n^A\otimes\mathcal{H}_m^A}^{\mathcal{H}_{n+m}^A}(N\otimes M)$
as an $\mathcal{H}_{n+m}^A$-module.

\begin{prop}\label{prop:bkB}
Let $J_1$ and $J_2$ be $\Z\rtimes\Z_2$-invariant subsets of
$\C^*$ such that $J_1 \cap J_2=\emptyset$.
\begin{enumerate}[{\rm(1)}]
\item
Let $N$ be an irreducible $\mathcal{H}_n^B$-module
of type $J_1$ and
$M$ an irreducible $\mathcal{H}_m^B$-module of type $J_2$.
Then $N\diamond M$ is an irreducible $\mathcal{H}_{n+m}^B$-module
of type $J_1 \cup J_2$.
\item Conversely if $L$
is an irreducible $\mathcal{H}_n^B$-module
of type $J_1 \cup J_2$, then there exist an integer $m$ $(0\le m\le n)$,
an irreducible $\mathcal{H}_m^B$-module $N$
of type $J_1$ and
an irreducible $\mathcal{H}_{n-m}^B$-module $M$ of type $J_2$
such that $L\simeq N\diamond M$.
\item
Assume that a $\Z\rtimes\Z_2$-orbit $J$ decomposes into
$J=J_+\sqcup J_-$ where $J_\pm$ are $\Z$-orbits and $J_-=(J_+)^{-1}$.
Assume that $\pm1,\pm p_0\not\in J$.
Then for any irreducible $\mathcal{H}_n^B$-module $L$ of type $J$,
there exists an irreducible $\mathcal{H}_n^A$-module $M$ such that
$L\simeq\Ind_{\mathcal{H}_n^A}^{\mathcal{H}_n^B}M$.
\end{enumerate}
\end{prop}
\begin{proof}
(1) \ Let $(N \diamond M)_{J_1,J_2}$ be the generalized eigenspace, 
where the eigenvalues of $X_i \ (1 \le i \le n)$ are in $J_1$ 
and the eigenvalues of $X_j \ (n<j \le n+m)$ are in $J_2$. 
Then $(N \diamond M)_{J_1,J_2}=N \otimes M$ by $J_1 \cap J_2=\emptyset$ by the above lemma and the shuffle lemma (e.g.\ \cite[Lemma 5.5]{G}). 
Suppose there exists non-zero $\mathcal{H}_{n+m}^B$-submodule $S$ 
in $N \diamond M$. 
Then $S_{J_1,J_2} \neq 0$ as an 
$\mathcal{H}_n^B \otimes \mathcal{H}_m^B$-module. 
Hence $S_{J_1,J_2}=N \otimes M$ by the irreducibility of $N \otimes M$ 
as an $\mathcal{H}_n^B \otimes \mathcal{H}_m^B$-module. 
We obtain $S=N \diamond M$.

\noindent
(2) \ For an irreducible $\mathcal{H}_{n}^B$-module $L$, the $\mathcal{H}_m^B \otimes \mathcal{H}_{n-m}^B$-module $L_{J_1,J_2}$ does not vanish for some $m$. Take an irreducible $\mathcal{H}_m^B \otimes \mathcal{H}_{n-m}^B$-submodule $S$ in $L$. Then there exist an irreducible $\mathcal{H}_m^B$-module $N$
of type $J_1$ and an irreducible $\mathcal{H}_{n-m}^B$-module $M$ of type $J_2$ such that $S=N \otimes M$. Hence there exists a surjective homomorphism $\Ind(N \otimes M)=N \diamond M \to L$. Since $N \diamond M$ is irreducible, this is an isomorphism. 

\noindent
(3) \ See \cite[Section 6]{M}.
\end{proof}
Hence in order to study $\mathcal{H}^B$-modules, it is enough to
study irreducible modules of type $J$ for  a $\Z\rtimes\Z_2$-orbit
$J$ in $\C^*$ such that $J$ is a $\Z$-orbit or $J$
contains one of $\pm1,\pm p_0$.

\subsection{The $a$-restriction and $a$-induction}
\begin{dfn}\label{def:Bef}
For $a \in \bb{C}^*$ and a finite-dimensional $\mathcal{H}_n^B$-module $M$, let us define the functors 
\begin{eqnarray*}
E_a:\mathcal{H}_n^B\fmod \to \mathcal{H}_{n-1}^B\fmod, \quad 
F_a:\mathcal{H}_n^B\fmod \to \mathcal{H}_{n+1}^B\fmod
\end{eqnarray*}
by\/{\rm:}
$E_aM$ is the generalized $a$-eigenspace of $M$ with respect 
to the action of $X_n$, and
$$F_aM\seteq
\Ind_{\mathcal{H}_n^B \otimes \bb{C}[X_{n+1}^{\pm{1}}]}%
^{\mathcal{H}_{n+1}^{B}}M \boxtimes \langle a \rangle,$$
where $\langle a \rangle$ is 
the $1$-dimensional representation of $\bb{C}[X_{n+1}^{\pm{1}}]$ 
defined by $X_{n+1} \mapsto a$. 

Define
\[
\widetilde{E}_aM\seteq\soc{E_aM}, \quad \widetilde{F}_aM\seteq\cosoc{F_aM}
\]
for $a \in \bb{C}^*$.
\end{dfn}
\begin{thm}[Miemietz \cite{M}]
Suppose $M$ is irreducible. Then $\widetilde{F}_aM$ is irreducible and $\widetilde{E}_aM$ is irreducible or $0$ for any $a \in \bb{C}^* \backslash 
\{\pm{1}\}$.
\end{thm}

\subsection{LLTA type conjectures for type B}
Now we take the case
\[
J=\set{p_1^k}{ k\in \bb{Z}_{\text{odd}}}.
\]
Assume that any of $\pm{1}$ and $\pm{p_0}$ is not contained in $J$. 
For short, we shall write $E_i$, $\widetilde{E}_i$, $F_i$ 
and $\widetilde{F}_i$ for $E_{p^i}$, $\widetilde{E}_{p^i}$, $F_{p^i}$ 
and $\widetilde{F}_{p^i}$, respectively. 

\begin{conj} 
\be[{\rm(1)}]
\item There are complete representatives 
\[
\set{L_b}{b \in B_\theta(0)}
\]
of the finite-dimensional irreducible $\mathcal{H}^B$-modules
of type $J$ such that 
\[
\widetilde{E}_{i}L_b=L_{\widetilde{E}_ib}, \quad \widetilde{F}_{i}L_b=L_{\widetilde{F}_ib}
\]
for any $i \in I\seteq\bb{Z}_{\text{odd}}$. 
\item For any $i \in \bb{Z}_{\text{odd}}$, 
let us define $E_{i,b,b'}, F_{i,b,b'} \in \bb{C}[q,q^{-1}]$ 
by the coefficients of the following expansions:
\begin{eqnarray*}
E_i\,G_\theta^{\text{up}}(b)=\sum_{b' \in B_\theta(0)}E_{i,b,b'}G_\theta^{\text{up}}(b'), \quad 
F_i\,G_\theta^{\text{up}}(b)=\sum_{b' \in B_\theta(0)}F_{i,b,b'}G_\theta^{\text{up}}(b'). 
\end{eqnarray*}
Then
\[
[E_iL_b:L_{b'}]=E_{i,b,b'}\vert_{q=1}, 
\quad [F_iL_b:L_{b'}]=F_{i,b,b'}\vert_{q=1}.
\]
Here $[M:N]$ is the composition multiplicity of $N$ in $M$ on $K_J^B$.
\ee
\end{conj}

\begin{rem}
There is a one-to-one correspondence between 
the above index set $B_\theta(0)$ 
and Syu Kato's parametrization (\cite{Kat}) 
of irreducible representations of $\mathcal{H}_n^B$ of type $J$. 
\end{rem}

\begin{rem}
\bnum
\item
For conjectures for other $\Z\rtimes\bb{Z}_2$-orbits $J$,
see \cite{EK}.
\item
Similar conjectures for type D are presented by
the second author and Vanessa Miemietz (\cite{KM}).
\ee
\end{rem}

\vspace{3ex}
\noindent
{\em Errata to ``Symmetric crystals and affine Hecke algebras of type B,}
Proc.\ Japan Acad., 82, no.\ 8, 2006, 131--136''
:
\bnum
\item
In Conjecture 3.8, $\lambda=\Lambda_{p_0}+\Lambda_{p_0^{-1}}$ should be 
read as 
$\lambda=\ssum_{a\in A}\Lambda_{a}$, where
$A=I\cap\{p_0,p_0^{-1},-p_0,-p_0^{-1}\}$.
We thank S.\ Ariki who informed us that
the original conjecture is false.
\item
In the two diagrams of $\Bt[\la]$ at the end of \S\,2,
$\lambda$ should be $0$.
\item
Throughout the paper, $A^{(1)}_\ell$ should be read as $A^{(1)}_{\ell-1}$.
\enum

\end{document}